\documentstyle[amssymb,amsfonts]{amsart}

\def\R{{\hbox{\bf R}}}
\def\H{{\hbox{\bf H}}}
\def\one{{\hbox{\rm \bf 1}}}
\def\osc{{\hbox{osc}}}

\def\J{{\cal J}}

\def\Z{{\hbox{\bf Z}}}
\def\eps{\varepsilon}
\def\pp{{p^\prime}}

\def\kp{{k^\prime}}
\def\jp{{j^\prime}}

\def\psI{{\psi_{2^s J}}}

\def\emph#1{{\it #1}}
\def\textbf#1{{\bf #1}}
\newenvironment{proof}{\noindent {\bf Proof} }{\endprf\par}
\def \endprf{\hfill  {\vrule height6pt width6pt depth0pt}\medskip}
\def\be#1{\begin{equation} \label{#1}}
\def\bs{\begin{split}}
\def\ba{\begin{align}}
\def\bas{\begin{align*}}
\def\eas{\end{align*}}

\theoremstyle{plain}
  \newtheorem{theorem}[subsection]{Theorem}
  \newtheorem{proposition}[subsection]{Proposition}
  \newtheorem{lemma}[subsection]{Lemma}
  \newtheorem{corollary}[subsection]{Corollary}

\theoremstyle{remark}

\theoremstyle{definition}
  \newtheorem{definition}[subsection]{Definition}

\begin{document}

\title[Homogeneous convolution operators]{The weak-type (1,1) of $L \log L$
homogeneous convolution operators}
\author{Terence Tao}
\address{Department of Mathematics, UCLA, Los Angeles, CA 90024}
\email{tao@@math.ucla.edu}

\begin{abstract}
We show that a homogeneous convolution kernel on an arbitrary homogeneous
group which is $L \log L$ on the unit annulus is bounded on $L^p$ for
$1 < p < \infty$ and is of weak-type $(1,1)$, generalizing the result of 
Seeger \cite{seeger:rough}.  The proof is in a similar spirit to that of
Christ and Rubio de Francia \cite{christ:rough}.
\end{abstract}

\maketitle

\section{Introduction}

Let $K$ be a homogeneous convolution
kernel on a homogeneous group $\H$, so that
so that
$$ K(t \circ x) = t^{-N} K(x)$$
for all $x \in \H$, $t > 0$, where $t \circ x$ is the dilation operation on $H$ and $t^N$ is the Jacobian of $x \mapsto t \circ x$.  Let $K_0$ be the restriction of $K$ to the unit annulus $A_0
= \{x \in \H: 1 \leq \rho(x) \leq 2 \}$, where $\rho$ is a norm associated
to the dilation structure.

We consider the question of what the minimal conditions on $K_0$ are
so that the convolution operator $T: f \mapsto f * K$ is
of weak-type $(1,1)$.  Here of course
$$ f*K(x) = \int f(y) K(y^{-1} x)\ dy.$$

Among the necessary conditions known are that $K_0$ must be in
$L^1(A_0)$, and $T$ must be bounded on $L^2$.  This in turn
necessitates that $K_0$ must have mean zero.
When $\H$ is an isotropic Euclidean space, the classical theorem
of Calder\'on and Zygmund \cite{calderon:rotations} shows that $T$ is
indeed bounded on $L^2$ when
$K_0$ has mean zero and is either odd and in
$L^1(A_0)$, or even and in the Orlicz space $L \log L(A_0)$.
In particular, we have boundedness on
$L^2$ whenever $K_0$ has mean zero and in $L \log L$.  This last condition
has been relaxed to $H^1$ and beyond; see \cite{stefanov:l2}.
Unfortunately, these arguments rely on the method of rotations and therefore
cannot be applied directly to the question of weak (1,1) boundedness
(cf. the discussion by R. Fefferman\cite{rfeff:entropy}).

It is natural to conjecture that analogous $L^2$ results hold for
arbitrary homogeneous groups.  If $K_0$ is odd and in $L^1$ one can use
the method of rotations and the work of Ricci and Stein
\cite{ricci:nilpotent} on convolution
operators on singular sets in homogeneous groups to obtain $L^2$
boundedness.  In case when $K_0$ is merely in $L \log L$ 
we shall use a variant of Littlewood-Paley theory
and an iterated $TT^*$ method to answer this conjecture affirmatively:

\begin{theorem}\label{L2}  If $K_0$ is in $L \log L$ and has mean zero,
then $T$ is bounded on $L^2$.
\end{theorem}

The weak-type $(1,1)$ question in Euclidean space has been considered
by several authors (\cite{christ:weak-1},\cite{christ:rough},
\cite{hofmann:weak},\cite{seeger:rough}); recently A. Seeger
\cite{seeger:rough} has shown that
$T$ is of weak-type $(1,1)$ on Euclidean space
whenever $K_0$ is in $L \log L$ and has mean 
zero.  The
corresponding questions for odd $L^1$ or even $H^1$ kernels remain
open.  We remark that the corresponding $(H^1,L^1)$ conjecture
is false by an example of Mike Christ.  For this and further discussion, 
see the survey in \cite{stefanov:l2}.

In this paper we generalize the result in \cite{seeger:rough} to
arbitrary homogeneous groups.  Specifically, we show that

\begin{theorem}\label{weak-11}  If $K_0$ is in $L \log L$ and $T$ is bounded
on $L^2$, then $T$ is of weak-type $(1,1)$.
\end{theorem}

Combining the two results and using duality we thus have

\begin{corollary}  If $K_0$ is in $L \log L$ and has mean zero, then
$T$ is bounded on $L^p$ for all $1 < p < \infty$, and is of weak-type
$(1,1)$.
\end{corollary}

The methods in \cite{seeger:rough} rely on the Euclidean Fourier transform and do not
appear to be adaptable to non-abelian settings.  Our approach is more in
the spirit of
Christ and Rubio de Francia \cite{christ:rough}, in that one considers
the expression $T f$ as an operator acting on the kernel $K_0$
rather than one acting on $f$.  One can then reduce weak $(1,1)$
boundedness to something resembling a $(L^2,L^2)$ estimate, which is now
amenable to orthogonality techniques such as the $TT^*$ and $(TT^*)^M$
methods.

The argument can also be used to treat the slightly smoother maximal
and square function operators corresponding to $L \log L$ generators
$K_0$, either by direct modification of the proof, or by using a
Radamacher function argument based on the fact that $\sum_i r_i(t) f * K_i$
is of weak-type $L^1$ uniformly in $t$.  For these operators the mean zero
condition is not required.  One can also use the arguments
in \cite{seeger:rough} to weaken the radial regularity on $K$ to a Dini-type
continuity condition.  We will not pursue these matters here.

This work originated from the author's dissertation at Princeton University under the inspiring guidance of Eli Stein.  The author is supported by NSF grant 9706764.

\section{Notation}

We will work exclusively with real-valued functions; none of our functions
will be complex-valued.

The letters $C$, (resp. $c$, $\epsilon$) will always be used to denote 
large (resp. small) positive constants that depend only on the 
homogeneous group $\H$ and any other specified quantities.  The values of these 
constants will change from line to line.  We use $A \lesssim B$ to denote
the statement that $A \leq CB$, and $A \sim B$ to denote the statement that
$A \lesssim B$ and $B \lesssim A$.

We define a 
homogeneous group to be a nilpotent Lie group $\H = \R^n$ 
with multiplication, inverse, dilation, and norm structures
$$ (x,y) \mapsto xy, \quad x \mapsto x^{-1}, \quad
(t,x) \mapsto t \circ x, \quad x \mapsto \rho(x)$$
for $x,y \in \H$, $t>0$,
where the multiplication and inverse operations are polynomial and form
a group with identity $0$, the dilation structure preserves the group
operations and is given in co-ordinates by
\be{dilate} t \circ (x_1, \ldots, x_n) = (t^{\alpha_1} x_1, \ldots, t^{\alpha_n} x_n)
\end{equation}
for some constants $0 < \alpha_1 \leq \alpha_2 \leq \ldots \leq \alpha_n$, and
$\rho(x)$ equals one on the Euclidean unit sphere, and satisfies
$\rho(t \circ x) = t \rho(x)$.  It can be shown that Lebesgue measure $dx$ is
a Haar measure for this group, and that $\rho(x) \sim \rho(x^{-1})$.  For
further properties of homogeneous groups see e.g. \cite{stein:large}.

We call $n$ the Euclidean dimension of $\H$, and the quantity 
$N = \alpha_1 + \ldots + \alpha_n$ the homogeneous dimension of 
$\H$.

We will always assume $\H$ to be a homogeneous group with Euclidean
dimension $n > 1$; the case $n=1$ can of course be treated by classical methods.
In addition to the homogeneous group structures mentioned above, we shall
also exploit the corresponding Euclidean structures
$$ (x,y) \mapsto x+y,\quad x \mapsto -x, \quad (t,x) \mapsto tx,\quad x \mapsto |x|,$$
together with the Euclidean inner product $(x,y) \mapsto x \cdot y$.

We shall also use the Euclidean structure of the exterior algebra $\Lambda$
of $\R^n$.  Recall that $\Lambda$ is spanned by basis elements of the form
$$ e_P = e_{p_1} \wedge \ldots \wedge e_{p_r}$$
where $0 \leq r \leq n$ and $P = (p_1, \ldots, p_r)$ is an increasing subsequence
of $1, \ldots, n$.  We give $\Lambda$ the usual inner product and norm structure
$$ (\sum_P a_P e_P) \cdot (\sum_Q b_Q e_Q) = \sum_P a_P b_P$$
and
$$ | \sum_P a_P e_P | = (\sum_P |a_P|^2)^{1/2}.$$
Later on we shall define some further structures on $\Lambda$ which are
more compatible with the non-isotropic dilation \eqref{dilate}.

We use $\prod_{i=1}^k x_i$ to denote the product
$x_1 \ldots x_k$, and $\prod_{i=k}^1 x_i$ to denote the product
$x_k \ldots x_1$.

We define a left-invariant
quasi-distance $d$ on $\H$ by $d(x,y) = \rho(x^{-1} y)$. 
A \emph{ball} $J = B(x_J, 2^j)$ with center $x_J$ and radius $2^j$ is 
defined to be any set of the form $$ J = \{ x: d(x, x_J) < 2^j\}$$
for some $x_J \in \H$ and $j \in \Z$.
If $J$ 
appears
in an expression, then $x_J$, $j$ are always understood to
be defined
as above.  If $C > 0$, then $CJ$ denotes the ball with the
same center as $J$ but $C$ times the radius.  
We use $J_\Delta$ to denote the annulus $CJ \backslash C^{-1} J$.

If $E$ is a finite set, we use $\# E$ to denote the cardinality of $E$; if $E$ is a 
measurable set, we use $|E|$ to denote the Lebesgue measure of $E$.  Note that 
$|t \circ E| = t^N |E|$ for all $t > 0$ and $E \subset \H$.

For each $t$
define the scaling map $\Delta[t]$ by
$$\Delta[t] f(y) = t^{-N} f(t^{-1} \circ y);$$
note that these operators are an isometry on $L^1$.

\section{Left-invariant differentiation structures}

Let $f(t)$ be a smooth function from $\R$ to $\H$.  The Euclidean derivative
$\partial_t f(t)$ can of course be defined by Newton's approximation
$$ f(t+\eps) = f(t) + \eps \partial_t f(t) + \eps^2 O(1)$$
for $\eps$ small.  We shall also need a left-invariant derivative $\partial^L_t f(t)$
defined by
$$ f(t+\eps) = f(t)(\eps \partial^L_t f(t)) + \eps^2 O(1).$$

If $f(t)$ is bounded, then the operation of left multiplication of $f(t)$ is
bilipschitz, and so we have
\be{comparable}
|\partial_t f(t)| \sim |\partial^L_t f(t)| \hbox{ whenever } |f(t)| \lesssim 1.
\end{equation}

We observe the product rule
\be{group-diff} \partial^L_t (f(t)g(t)) = \partial^L_t g(t) +
C[g(t)] \partial^L_t f(t)
\end{equation}
where the linear transformation $C[x]: \R^n \to \R^n$ is the derivative
of the conjugation map $y \mapsto x^{-1} y x$ at the origin.  In other
words, for all $x \in \H$, $v \in \R^n$ we have
$$ x^{-1}(\eps v)x = \eps C[x] v + \eps^2 O(1).$$
The rule \eqref{group-diff} is easily verified by expanding 
$f(t+\eps)g(t+\eps)$ to first order in two different ways.  We note the
identities
\be{c-twine} C[t \circ x] (t \circ v) = t \circ (C[x] v),
\quad C[x]^{-1} = C[x^{-1}].
\end{equation}
Since $C[x]$ and its inverse are both polynomial in $x$, we have
\be{c-bound}
|C[x] v| \sim |v| \hbox{ whenever } |x| \lesssim 1.
\end{equation}

Now suppose $F(x)$ is a smooth function from $\R^n$ to $\H$.  We define the
left-invariant derivative $D^L_x F(x)$ to be the matrix with columns given by
$$D^L_x F(x) = (\partial^L_{x_1} F(x), \ldots, \partial^L_{x_n} F(x)).$$
In other words, we have the Newton approximation
$$ F(x + \eps v) = F(x) (\eps D^L_x F(x) v) + \eps^2 O(1).$$
Since $dx$ is a Haar measure, we see that the determinant of $D^L_x F(x)$
is equal to the Jacobian of $F$ at $x$ with respect to Lebesgue measure.
We note that
\be{det-form}
\det D^L_x F(x) (e_1 \wedge \ldots \wedge e_n)
= \partial^L_{x_1} F(x) \wedge \ldots \wedge \partial^L_{x_n} F(x).
\end{equation}

The vector field
$$X(x) = \partial^L_t (t \circ x)|_{t = 1}.$$
shall be crucial in our arguments.  An equivalent definition is
\be{x-def} 
(1+\eps) \circ x = x (\eps X(x)) + \eps^2 O(1).
\end{equation}

Note that $X$ commutes with dilation:
\be{x-dil} X(t \circ x) = t \circ X(x).
\end{equation}
Since $X$ depends polynomially on $x$, we therefore have
\be{x-bound}
\rho(X(x)) \lesssim \rho(x)
\end{equation}
For comparison, we also observe the bound
\be{trivx-bound}
\rho(\partial_t (t \circ x)) \sim \rho(t \circ x)
\end{equation}
which follows immediately from \eqref{dilate}.

The left-invariant derivative $\partial^L_t$ interacts with dilations via
the formula
\be{dil-diff} 
\partial^L_t (s(t) \circ f(t))
= s(t) \circ \partial^L_t f(t)
+ \frac{s'(t)}{s(t)} (s(t) \circ X[f(t)])
\end{equation}
which is verified by expanding $s(t+\eps) \circ f(t+\eps)$ to first
order in two different ways.  Finally, we have

\begin{lemma}\label{x-invert} 
The map $X: \R^n \to \R^n$ is a polynomial diffeomorphism with
Jacobian comparable to 1.
\end{lemma}

\begin{proof}  
From the monotonicity assumptions on $\alpha_i$ and the assumption that
dilations preserve the multiplication structure, it is easy to see
that the multiplication law $(x,y) \mapsto xy$ on $\H$ must have
the upper diagonal form
\be{mult}
\begin{split}
 (xy)_1 &= x_1 + y_1\\
 (xy)_2 &= x_2 + y_2 + P_2(x_1,y_1)\\
 (xy)_3 &= x_3 + y_3 + P_3(x_1,x_2,y_1,y_2)\\
&\ldots \\
 (xy)_n &= x_n + y_n + P_n(x_1, \ldots, x_{n-1}, y_1, \ldots, y_{n-1})
\end{split}
\end{equation}
where $P_2, \ldots, P_n$ are polynomials.

From \eqref{dilate} and \eqref{x-def} we have
$$ \eps(\alpha_1 x_1, \ldots, \alpha_n x_n) = x (\eps X(x)) + \eps^2 O(1).$$
Inserting this into \eqref{mult} and solving recursively for the components
of $X(x)$ we see that
\bas
X(x)_1 &= \alpha_1 x_1\\
X(x)_2 &= \alpha_2 x_2 + Q_2(x_1)\\
&\ldots\\
X(x)_n &= \alpha_n x_n + Q_n(x_1, \ldots, x_{n-1})
\end{align*}
for some polynomials $Q_2, \ldots, Q_n$ which depend on the
$\alpha_i$.  The claim follows.
\end{proof}

\section{Proof of Theorem \ref{L2}.  Kernel truncation and frequency
localization.}

We now begin the proof of Theorem \ref{L2}.  The heart of the
argument is an iterated $TT^*$ method, in the spirit
of Christ and Rubio de Francia \cite{christ:rough}.
We may normalize
$\|K_0\|_{L \log L} = 1$.

We first partition the kernel $K$ dyadically.  From the identity
$$ K = \frac{1}{\ln 2} \int \Delta[t] K_0\ \frac{dt}{t}$$
we have the decomposition $K = \sum_j S_j K_0$, where
$S_j$ is the operator
\begin{equation}\label{kj-def}
S_j F = 2^{-j} \int \varphi(2^{-j} t) \Delta[t] F\ dt,
\end{equation}
and $\varphi$ is a bump function adapted to $\{t \sim 1\}$ such that
$\sum_j 2^{-j} t \varphi(2^{-j} t) = \frac{1}{\ln 2}$.  Note that 
\be{sj-est}
\|S_j F\|_1 \lesssim \|F\|_1
\end{equation}
uniformly in $j$.  From the a priori assumptions on $K_0$ we see that
the $S_j K_0$ are all $C^\infty_0$ functions.

We need to show that
$$ \| f * \sum_j S_j K_0 \|_2 \lesssim \|f\|_2.$$

Write $K_0 = \sum_{s \geq 0} K_0^s$, where 
$$K_0^s = \hat K_0^s - \frac{\chi_{A_0}}{|A_0|} \int_{A_0} \hat K_0^s$$
and $\hat K_0^s$ is the portion of $K_0$ on the set 
$$2^{2^s} \leq 1+|K_0| < 2^{2^{s+1}}.$$
Note that each $K_0^s$ has mean zero.
By the triangle inequality and the computation
$$ \sum_{s \geq 0} 2^s \|K_0^s\|_1 \lesssim 
\sum_{s \geq 0} 2^s \|\hat K_0^s\|_1 \lesssim 
\|K_0\|_{L \log L} = 1$$
it thus suffices to show that
$$ \| f * \sum_j S_j K_0^s \|_2 \lesssim \|f\|_2 (2^s \|K_0^s\|_1 + 
2^s 2^{-\eps 2^s})$$
for all $s \geq 0$.

Fix $s$.  For each integer $k$, let $T_k$ denote the operator
$$ T_k f = f * \sum_{j=k2^s}^{(k+1)2^s - 1} S_j K_0^s.$$
Our task is to show the operator norm estimate
$$ \| \sum_k T_k \| \lesssim 2^s \|K_0^s\|_1 + 2^s 2^{-\eps 2^s}.$$
From Young's inequality and \eqref{sj-est} we  have
$$ \| T_k f \|_2 \lesssim \|f\|_2 2^s \|K_0^s\|_1$$
for all integers $k$.
In particular, we have the operator norm estimates
$$ \| T_k T_\kp^* \|, \| T_k^* T_\kp \|
\lesssim (2^s \|K_0^s\|_1)^2$$
for all $k,\kp$.  We will shortly show that
\be{L2-targ}
\| T_k T_{\kp}^* \|, \| T_k^* T_\kp \|
\lesssim 2^{2s} 2^{-\eps 2^s |k-\kp|}
\end{equation}
for $|k-\kp| \geq C$, where $C$ is a large constant to be
determined later.  From these estimates the desired
bound on $\| \sum_k T_k \|$ follows from the Cotlar-Knapp-Stein
lemma (see e.g. \cite{stein:large}).

It remains to prove \eqref{L2-targ}.  We prove only the first estimate,
as the second is analogous.  We rewrite this as

$$ \|\sum_{j=k2^s}^{(k+1)2^s-1} \sum_{j'=k'2^s}^{(k'+1)2^s-1}
f * S_j K_0^s * S_{j'} \tilde K_0^s \|_2 \lesssim 2^{2s} 2^{-\eps 2^s |k-\kp|} \|f\|_2$$
where $\tilde F$ denotes the function $\tilde F(x) = F(x^{-1})$.  

By the triangle inequality it suffices to show that
\be{decay}
 \| f * S_j K_0^s * S_{j'} \tilde K_0^s \|_2 \lesssim 2^{-\eps |j-\jp|}
\|f\|_2
\end{equation}
for all integers $j,\jp$ for which $|j-\jp| > C2^s$.  

The next step is to introduce a form of Littlewood-Paley theory, although
we shall avoid any explicit use of the Fourier transform.
Fix a function $\phi$ on the unit ball with $\|\phi\|_{C^1} 
\lesssim 1$ which has unit mass.  We may also assume that 
$\phi = \tilde \phi$.  For each integer $k$, write
$$ \Psi_k = \Delta[2^{k-1}] \phi - \Delta[2^k] \phi.$$
Note that $\Psi_k$ is supported on the ball of radius $C2^{k}$, has mean zero,
and $\tilde \Psi_k = \Psi_k$.

Since $1 = \sum_k \Psi_k$, we may write 
$$ f * S_j K_0^s * S_0 \tilde K_0^s 
= \sum_k \sum_\kp f * S_j K_0^s * \Psi_k * \Psi_\kp * S_{\jp} \tilde K_0^s.$$
Suppose for the moment that we could prove

\begin{proposition}\label{L2-key}  For any integers $j$, $k$, and
any $L^\infty$ function $K_0$ on the unit annulus with mean zero, we have
$$ \| f * S_j K_0 * \Psi_k \|_2 \lesssim 2^{-\eps |j-k|} \|f\|_2 
\|K_0\|_\infty.$$
\end{proposition}

Then we would have the estimates
$$ \| f * S_j K_0^s * \Psi_k \|_2 \lesssim 2^{2^{s+1}} 2^{-\eps |j-k|} \|g\|_2.$$
and (by duality)
$$ \| g * \Psi_\kp * S_\jp \tilde K_0^s\|_2 \lesssim 2^{2^{s+1}} 
2^{-\eps |\kp-\jp|} \|g\|_2.$$
Combining these two estimates we see that
\be{first-fss}
 \|  f * S_j K_0^s * \Psi_k * \Psi_\kp * S_\jp \tilde K_0^s \|_2
\lesssim 2^{2^{s+2}} 2^{-\eps |j-k|} 2^{-\eps |\kp-\jp|} \|f\|_2.
\end{equation}
On the other hand, by Young's inequality we also have the estimate
$$ \|  f * S_j K_0^s * \Psi_k * \Psi_\kp * S_\jp \tilde K_0^s \|_2
\lesssim \|f\|_2 \|S_j K_0^s\|_1 \|\Psi_k * \Psi_\kp\|_1
\|S_0 \tilde K_0^s\|_1.$$
Since $K_0$ is in $L \log L$, $K_0^s$ is in $L^1$, and so by \eqref{sj-est} we have
$$ \|S_j K_0^s\|_1, \|S_\jp \tilde K_0^s\|_1 \lesssim 1.$$
Also, from the smoothness and mean zero conditions on $\Psi_k$, $\Psi_\kp$
we have
$$ \|\Psi_k * \Psi_\kp\|_1 \lesssim 2^{-\eps |k-\kp|}.$$
Thus we obtain the bound of
$$ \|  f * S_j K_0^s * \Psi_k * \Psi_\kp * S_\jp \tilde K_0^s \|_2
\lesssim 2^{-\eps |k-\kp|} \|f\|_2.$$
Taking the geometric mean of this with \eqref{first-fss} we obtain
$$  \|  f * S_j K_0^s * \Psi_k * \Psi_\kp * S_\jp \tilde K_0^s \|_2
\lesssim 2^{2^{s+1}} 2^{-\eps |j-k|/2} 2^{-\eps|k-\kp|/2} 2^{-\eps |\kp-\jp|/2} 
\|f\|_2.$$
If we then sum this in $k$ and $\kp$ we obtain 
$$ \| f * S_j K_0^s * S_\jp \tilde K_0^s \|_2
\lesssim 2^{2^{s+1}} |j-\jp| 2^{-\eps |j-\jp|/2},$$
which gives \eqref{decay} for some $\eps > 0$,
if $|j-\jp| > C2^s$ for a sufficiently large $C$.

\section{Proof of Theorem \ref{L2} continued.  Iterated $TT^*$ methods.}

It thus remains to prove Proposition \ref{L2-key}.  We may normalize
so that $\|K_0\|_\infty = 1$; by scale invariance we may assume that
$j=0$.  

If $k \geq -C$ then from the mean zero condition on $K_0$ and
the smoothness of $\Psi_k$ we have
$$ \| S_0 K_0 * \Psi_k \|_1 \lesssim 2^{-\eps k}$$
and the desired bound thus follows from Young's inequality.
We may therefore assume that $k < -C$.

Fix $k = -s$ for some $s>C$.  Our task is now to show
$$ \| f * S_0 K_0 * \Psi_{-s} \|_2 \lesssim 2^{-\eps s} \|f\|_2.$$

It is possible to use Fourier techniques to handle this estimate, taking advantage of the microlocal
regularity properties of $S_0 K_0$.  However, we shall
pursue a different approach based on the iterated $T^*T$ method,
as we shall need these techniques later on for the (more difficult)
weak (1,1) estimate.

Roughly speaking, the idea is as follows.  The kernel $S_0 K_0$
is smooth along the ``radial'' direction, but is otherwise rough.
Thus there is no obvious way to exploit the cancellation properties
of $\Psi_{-s}$.  However, if we convolve $S_0 K_0$ with itself $n$ times then
one should obtain a kernel which is smooth in $n$ separate directions
at any given point.  Assuming that these directions are linearly
independent, the iterated kernel thus has isotropic regularity, and
one will pick up the desired $2^{\eps s}$ gain by exploiting the moment
conditions of $\Psi_{-s}$.  Of course, there will be an exceptional portion
of the convolution in which the directions of smoothing are not
independent.  For this portion one cannot exploit cancellation and
one must instead replace everything by absolute values.

We now turn to the details. By the $T^*T$ method, it suffices to show that
$$ \| f * \Psi_{-s} * S_0 \tilde K_0 * S_0 K_0 * \Psi_{-s} \|_2 \lesssim 
2^{-\eps s} \|f\|_2.$$
From the operator norm identity $\|T^* T\| = \|(T^* T)^n\|^{1/n}$,
it thus suffices to show that
$$ \| f * \Psi_{-s} * S_0 \tilde K_0 * S_0 K_0 * \Psi_{-s}
* \ldots * \Psi_{-s} * S_0 \tilde K_0 * S_0 K_0 * \Psi_{-s} \|_2
\lesssim 2^{-\eps s} \|f\|_2$$
for a slightly different value of $\eps > 0$, where the convolution
is iterated $n = \dim \H$ times.
By Young's inequality it suffices to show that
$$ \| \Psi_{-s} * S_0 \tilde K_0 * S_0 K_0 * \Psi_{-s}
* \ldots * \Psi_{-s} * S_0 \tilde K_0 * S_0 K_0 * \Psi_{-s} \|_1
\lesssim 2^{-\eps s}.$$

The function $\Psi_{-s} * S_0 \tilde K_0$ is bounded in $L^1$,
and is therefore an average of delta functions in $B(0,C)$.  From
Minkowski's inequality, it
therefore suffices to show that
$$ \| \delta_{w_1} * S_0 K_0 * \delta_{w_2} * S_0 K_0
* \ldots * \delta_{w_n} * S_0 K_0 * \Psi_{-s} \|_1
\lesssim 2^{-\eps s}$$
uniformly for all $w_1, \ldots, w_n \in B(0,C)$.

Fix $w=(w_1, \ldots, w_n)$.
It suffices
to show that
$$ |\langle \delta_{w_1} * S_0 K_0 * \delta_{w_2} * S_0 K_0
* \ldots * \delta_{w_n} * S_0 K_0 * \Psi_{-s},g \rangle| \lesssim 2^{-\eps s}$$
for all test functions $g$ which are normalized in $L^\infty$.

Fix $g$.  We write the left-hand side as
$$ |\int \int \int \Psi_{-s}(x)
g(\Phi_y(t) x)
\prod_{q=1}^n K_0(y_q) \varphi(t_q)\ dy dt dx|$$
where $t = (t_1, \ldots, t_n) \in [C^{-1},C]^n$, $y = (y_1, \ldots, y_n) \in A_0^n$, and
\be{Phi-def}
 \Phi_{y}(t) = \prod_{q=1}^n w_q (t_q \circ y_q).
\end{equation}

The treatment of this integral depends on whether the map $\Phi_y: \R^n \to \H$
is degenerate or not.  This degeneracy is measured by the Jacobian $\det D^L_t(\Phi_y(t))$.  Accordingly, we split our estimates into
\be{non-cancel}
|\int \int \int \Psi_{-s}(x) \eta(2^{n\eps s} \det D^L_t(\Phi_{y})(t))
g(\Phi_{y}(t)x)
\prod_{q=1}^n K_0(y_q) \varphi(t_q)\ dy dt dx|
\lesssim 2^{-\eps s}
\end{equation}
and
\be{cancel}
|\int \int \int \Psi_{-s}(x) [1-\eta(2^{n\eps s} \det D^L_t(\Phi_{y})(t))]
g(\Phi_{y}(t)x)
\prod_{q=1}^n K_0(y_q) \varphi(t_q)\ dy dt dx|
\lesssim 2^{-\eps s},
\end{equation}
where $\eta$ is a smooth non-negative bump function which equals $1$ near 1.

\section{Proof of Theorem \ref{L2} continued.  The degenerate portion
of the integral.}

To show \eqref{non-cancel} we simply replace everything by
absolute values, and use the bounds on 
$K_0$, $g$, and $\varphi$ to reduce to
$$ \int \int_{[C^{-1},C]^n} \int_{A_0^N} 
|\Psi_{-s}(x)| \eta(2^{n\eps s} \det D^L_t(\Phi_{y})(t))\ dy dt dx
\lesssim 2^{-\eps s}.$$
Performing the $x$ integration and taking supremums in the $t$ integral,
we reduce to
$$
\sup_{t \in [C^{-1},C]^n}
\int_{y \in A_0^n: |\det D^L_t(\Phi_{y})(t)| \lesssim 2^{-n\eps s}}
dy \lesssim 2^{-\eps s}.$$
Fix $t \in [C^{-1},C]^n$.  From \eqref{det-form} we have
$$ |\det D^L_t(\Phi_y)| = |\partial^L_{t_1} \Phi_y \wedge \ldots
\wedge \partial^L_{t_n} \Phi_y|.$$
From \eqref{Phi-def}, \eqref{group-diff} and \eqref{dil-diff} we have
\be{Phi-diff}
 \partial^L_{t_q} \Phi_y = t_q^{-1} C[Q_q] (t_q \circ X(y_q))
\end{equation}
for all $q$, where $Q_q$ is the quantity
$$ Q_q = \prod_{j=q+1}^n w_j (t_j \circ y_j).$$
Since $|Q_q| \lesssim 1$, $t_q \sim 1$, and $|y_q| \sim 1$, we see
from \eqref{c-bound} that 
$$
|\partial^L_{t_q} \Phi_y| \sim 1.
$$
We therefore have
\be{ei-decomp}
\int_{y \in A_0^n: |\det D_t(\Phi_{y})(t)| \lesssim 2^{-n\eps s}}
dy_1 \ldots dy_n \leq \sum_{q=1}^{n-1} |E_q|,
\end{equation}
where $E_q$ is the set.
$$ E_q = \{ y \in A_0^n:
| \partial^L_{t_q} \Phi_y \wedge \ldots
\wedge \partial^L_{t_n} \Phi_y |
< 2^{-\eps s}
| \partial^L_{t_{q+1}} \Phi_y \wedge \ldots
\wedge \partial^L_{t_n} \Phi_y | \}.$$
It thus suffices to show that $|E_q| \lesssim 2^{-\eps s}$
for each $q$.

Fix $q$, and freeze all the $y_j$ variables except for $y_q$.
From \eqref{Phi-diff} and \eqref{c-bound}, we see
that in order for $y$ to be in $E_q$, $X(y_q)$ must live in a boundedly finite union of $2^{-\eps s}$-neighbourhoods of planes.  These planes depend
only on
$Q_q$ and $\partial_{t_{q+1}} \Phi_y \wedge \ldots \wedge \partial_{t_n} \Phi_y$,
and so are independent of $y_q$.  Since $y_q$ is bounded, 
we thus see from Lemma \ref{x-invert} that $y_q$ lives in a finite
union of $C2^{-\eps s}$-neighbourhoods of compact hypersurfaces.
In particular, the variable $y_q$ must
range in a set of measure $O(2^{-\eps s})$.  The desired bound on $E_q$ follows
by unfreezing the remaining $y$ variables.  This concludes the proof of 
\eqref{non-cancel}.

\section{Proof of Theorem \ref{L2} continued.  The non-degenerate portion
of the integral.}

To finish the proof of Theorem \ref{L2} we must show \eqref{cancel}.
Since $K_0$ is in $L^\infty(A_0)$, it is in $L^1$, and it suffices to show

$$
|\int \int \Psi_{-s}(x) g(\Phi_{y}(t)x) [1-\eta(2^{n\eps s} \det D^L_T(\Phi_{y})(t))]
\prod_{i=1}^n \varphi(t_i)\ dt dx|
\lesssim 2^{-\eps s}
$$
uniformly in $y$.  

Fix $y$.  We now utilize the moment conditions in the $\Psi_{-s}$ by
rewriting $\Psi_{-s}$ as a (Euclidean) divergence of a function which is small in $L^1$.
More precisely, we shall use

\begin{lemma}\label{integ}  Let $f$ be a function on $B(0,C)$ with mean zero
and $\|f\|_1 \lesssim 1$.  Then there exists functions $f_1, \ldots, f_n$
supported on a slightly larger ball $B(0,C)$ with $\|f_i\|_1 \lesssim 1$ and
$$ f(x) = \sum_i \partial_{x_i} f_i(x).$$
\end{lemma}

\begin{proof}  Without loss of generality we may assume that $f$ is supported on the unit cube $[0,1]^n$.

When $n=1$ the lemma is clear.  For $n>1$ we write
$$ f(x_1, \ldots x_n) = \partial_{x_n} f_n(x_n) 
+ F_{x_n}(x_1, \ldots, x_{n-1})$$
for all $x \in [0,1]^n$, where 
$$ f_n(x_n) = \int_{x'_n \leq x_n} f(x')\ dx'$$
and 
$$ F_{x_n}(x_1, \ldots, x_{n-1}) = f(x_1, \ldots, x_n) - \int_{x'_n = x_n} f(x')\ dx.$$
Clearly $f_n$, $F_{x_n}$ have bounded $L^1$ norm on the unit cube, and $F_{x_n}$ has mean zero for each $x_n$.  The lemma then follows from induction.
\end{proof}

Applying this lemma to $f = \Psi_0$ and then rescaling, we may write
$$ \Psi_{-s}(x) = \sum_{i=1}^n \partial_{x_i} f_i(x)$$
where the functions $f_i$ are supported on $B(0,C 2^{-s})$ and satisfy
\be{l1-gain}
\|f_i\|_1 \lesssim 2^{-\alpha_i s}.
\end{equation}
  We thus need to show that
\be{integ-L2} | \int \int \partial_{x_i} f_i(x) 
g(\Phi_{y}(t)x) a(t)\ dt dx| \lesssim 2^{-\eps s}
\end{equation}
for all $i=1,2,\ldots, n$, where
$$ a(t) = [1-\eta(2^{n\eps s} \det D^L_T(\Phi_{y})(t))] \prod_{i=1}^n \varphi(t_i).$$

Fix $i$.  The idea is to use integration by parts to somehow move the 
derivative $\partial_{x_i}$ onto the smooth function $a$, so that one can
exploit \eqref{l1-gain} and the $L^\infty$ control on $g$.

If we integrate by parts in the $x_i$ variable, the left-hand side
of \eqref{integ-L2} becomes
$$
|  \int \int f_i(x) 
\partial_{x_i} g(\Phi_{y}(t)x) a(t)\ dt dx|.$$
From \eqref{l1-gain}, it thus suffices (if $\eps$ is chosen sufficiently small) to show that
\be{integ-l3}
| \int \partial_{x_i} g(\Phi_{y}(t)x) a(t)\ dt| 
\lesssim 2^{C \eps s},
\end{equation}
for all $x \in B(0,C)$.

Fix $x$.  We now apply the following application of the chain rule,
which allows one to convert a derivative of one variable to a derivative
on another variable, provided that a certain Jacobian is non-zero.

\begin{lemma}\label{chain}  Let $f: \R \times \R^n \to \H$ and $F: \H \to \R$ be
smooth functions.  Then
\be{chain-rule} \partial_{s} F(f(s,t)) = \nabla_t F(f(s,t)) \cdot  
 (D^L_t f(s,t))^{-1} \partial^L_{s} f(s,t)
\end{equation}
whenever $\det D^L_t f(s,t)$ is non-zero.
\end{lemma}

\begin{proof}  For any small $\eps$, we have the Newton approximations
$$ F(f(s+\eps,t)) = F(f(s,t)) + \eps \partial_s F(f(s,t)) + \eps^2 O(1)$$
$$ f(s+\eps,t) = f(s,t) (\eps \partial^L_s f(s,t)) + \eps^2 O(1)$$
$$ f(s,t+\eps v) = f(s,t) (\eps D^L_t f(s,t) v) + \eps^2 O(1)$$
$$ F(f(s,t + \eps v)) = F(f(s,t)) + \eps \nabla_t F(f(s,t)) \cdot v + \eps^2 O(1).$$
Combining all these estimates with $v = D^L_t f(x,t)^{-1} \partial^L_s f(s,t)$
and letting $\eps \to 0$ gives the result.
\end{proof}

From this lemma, \eqref{integ-l3} becomes
$$
| \int \nabla_t g(\Phi_y(t)x)
\cdot (D^L_t (\Phi_y(t) x))^{-1} \partial^L_{x_i}(\Phi_y(t)x) a(t)\ dt|
\lesssim 2^{C\eps s}.$$
By another integration by parts and the fact that $g \in L^\infty$, it 
suffices to show the uniform estimate
$$ | \nabla_t \cdot [(D^L_t (\Phi_y(t) x))^{-1} \partial^L_{x_i}(\Phi_y(t)x) a(t)] | \lesssim 2^{C\eps s}.$$
But this is easily verified, since all variables are compactly supported and all
functions are smooth, with norms at most $O(2^{C\eps s})$.  The 
$(D^L_t (\Phi_y(t) x))^{-1}$ term is well-behaved since
$$ |\det D^L_t(\Phi_y(t) x)| \sim |\det D^L_t(\Phi_y(t))| \gtrsim 2^{-n\eps s}$$
on the support of $a(t)$.
This completes the proof of Proposition
\ref{L2-key} and thus Theorem \ref{L2}.

\section{Proof of Theorem \ref{weak-11}.  Truncation of the kernel
and strong-type estimates.}

We now begin the proof of Theorem \ref{weak-11}.  The arguments
will be similar in flavor to the ones used to prove Theorem \ref{L2}, but with
two major differences.  Firstly, because the function $f$ is now
only controlled in $L^1$, one is forced (as in \cite{christ:rough})
to perform the $TT^*$ method with respect to $K_0$ rather than $f$.  
Secondly, the Littlewood-Paley operators are not particularly useful in the $L^1$
setting, and we cannot reduce to an estimate on a single scale 
such as Proposition \ref{L2-key}.  Instead, we are forced to consider
the interactions between several scales.  This will cause an increase
in complexity in our arguments.  We remark that if one were to treat
the maximal function or square function instead of the singular integral,
then one could again localize to a single scale; cf. the arguments
in \cite{christ:rough}.

Let $K$ be as in the statement of the theorem.
We wish to 
show that
$$ \{ |f * K| \gtrsim \alpha \} \lesssim \alpha^{-1} \|f\|_1 \|K_0\|_{L \log L}.$$
We may assume that $f$ is a $C^\infty_0$ function.  By linearity we may assume that
$\alpha = 1$ and $\|K_0\|_{L \log L} = 1$.  

We perform the standard Calder\'on-Zygmund decomposition of $f$ at height
$1$ to obtain $f = g + \sum_J b_J$, where $\|g\|_1 \lesssim \|f\|_1$, 
$\|g\|_\infty \lesssim 1$, the $J$ range over a collection of disjoint balls with $\sum_J |J| \lesssim \|f\|_1$, and for each $J$ the functions $b_J$
are supported on $CJ$ with
\be{bj-prop}
\|b_J\|_{L^1(CJ)} \lesssim |J|, \quad \int b_J = 0.
\end{equation}
Since $f$ is smooth, the collection of $J$ is finite.  We may arrange matters so that the $b_J$ are smooth.

We now proceed with the standard reduction argument as employed in
\cite{christ:rough}, \cite{seeger:rough}.  

As in the proof of Theorem \ref{L2}, we decompose $K = \sum_j S_j K_0$.
We need to estimate the set where
$$ f * K = g * K + \sum_{s \leq C} \sum_J b_J * S_{j+s} K_0 
+ \sum_{s > C}\sum_J b_J * S_{j+s} K_0 $$
is essentially greater than 1; here and in the sequel, $j = j(J)$ is the integer
such that $J$ has side-length $2^j$.

The first term can be handled by 
the $L^2$ boundedness hypothesis and Chebyshev's
inequality because $g$ is in $L^2$ with norm $O(\|f\|_1^{1/2})$.
The second term is supported in $\bigcup_J CJ$, and so
that contribution is acceptable since $\sum_J |CJ| \lesssim \|f\|_1$.

To handle the remaining term it suffices to show that
$$
|\{ \sum_{s > C} |\sum_J b_J * S_{j+s} K_0| \gtrsim 1\}|
\lesssim \sum_J |J|.
$$

We introduce a cutoff to emphasize the fact that $b_J * S_{j+s} K_0$
is supported on the annulus $(2^s J)_\Delta$.  Namely, we rewrite the above
as
\be{weak-est}
|\{ \sum_{s > C} |\sum_J \psi_J(b_J * S_{j+s} K_0)| \gtrsim 1\}|
\lesssim \sum_J |J|
\end{equation}
where $\psi_J(x) = \psi(2^{-j} \circ (x_J^{-1} x))$ and $\psi$ is 
a suitable cutoff function supported on (a slight thickening of)
the unit annulus $A_0$.  We now claim that \eqref{weak-est} will follow from

\begin{proposition}\label{first}  Let $s \geq C$, 
$\J$ be a non-empty finite collection
of disjoint balls such that 
\be{size}
\sum_J |J| \lesssim 1,
\end{equation}
and $b_J$ be a collection of smooth functions satisfying \eqref{bj-prop}.  Let $\psi_J$ be defined as above.  Let $1 < p < 2$ be an exponent.
Then there exists an exceptional set $E = E_s$ 
such that 
$|E| \lesssim 2^{-\eps s}$ and
\begin{equation}\label{first-reduction}
\| \sum_J \psi_J(b_J * S_{j+s} F_J) \|_{L^p(E^c)} \lesssim 2^{-\eps s} 
(\sum_J |J| \|F_J\|_2^2)^{1/2}
\end{equation}
for all functions $F_J$ in $L^2(\H)$.
\end{proposition}

We will prove this proposition in later sections.  For now,
we show why Proposition \ref{first} implies \eqref{weak-est}.  
It suffices by dilation invariance
to verify \eqref{weak-est} in the case when $\sum_J |J| \sim 1$.  In particular,
we may assume that \eqref{size} holds.

For each $s > C$ we decompose $K_0$ as
$K_0 = K^{\leq s} + K^{>s}$, where $K^{\leq s}$ is the portion of
$K_0$ supported on the set $|K_0| \lesssim 2^{\eps s/2}$.  We have to show
that
\be{former}
 |\{ |\sum_{s > C} \sum_J \psi_J(b_J * S_{j+s} K^{>s}) | \gtrsim 1\}| \lesssim 1
\end{equation}
and
\be{latter}
 |\{ |\sum_{s > C} \sum_J \psi_J(b_J * S_{j+s} K^{\leq s}) | \gtrsim 1\}| \lesssim 1.
\end{equation}
To show \eqref{former} it suffices by Chebyshev's inequality to show
that
$$\| \sum_{s > C} \sum_J \psi_J(b_J * S_{j+s} K^{> s}) \|_1 \lesssim 1.$$
But from \eqref{sj-est}, \eqref{bj-prop}, and Young's inequality one sees that
$$ \| \psi_J(b_J * S_{j+s} K^{> s}) \|_1 \lesssim |J| \|K^{> s}\|_1,$$
and the desired estimate follows from \eqref{size}
and the observation that 
$$\sum_{s > C} \|K^{> s}\|_1 \lesssim \|K_0\|_{L \log L} =  1.$$

To show \eqref{latter} it suffices by Chebyshev and the
observation
 $|\bigcup_{s > C} E_s| \lesssim 1$ to show that
$$ \| \sum_{s > C} \sum_J \psi_J(b_J * S_{j+s} K^{\leq s})\|_{L^p((\bigcup_{s>C} E_s)^c)} \lesssim 1.$$
By the triangle inequality it suffices to show
$$ \| \sum_J \psi_J(b_J * S_{j+s} K^{\leq s})\|_{L^p(E_s^c)} \lesssim 2^{-\eps s/2}$$
for each $s$.
But this follows from \eqref{first-reduction} with $F_J = K^{\leq s}$ for all $J$, since 
$$(\sum_J J \|K^{\leq s}\|_2^2)^{1/2} \lesssim 
\| K^{\leq s} \|_2 \lesssim \|K^{\leq s}\|_\infty \lesssim 2^{\eps s/2}.$$
This completes the derivation of \eqref{weak-est} from Proposition
\ref{first}.  It remains only to prove Proposition \ref{first}.

\section{Proof of Theorem \ref{weak-11} continued.  Bounded overlap of
dilated balls.}

In the remainder of the argument, $s>C$ and $1 < p < 2$ will be fixed.

To prove Proposition \ref{first}, we first prove it under a natural
multiplicity assumption on the overlap of the sets $2^s J$.  More precisely, we will show in later sections that

\begin{proposition}\label{next}  Let 
$\J$ be a non-empty finite collection
of disjoint balls such that \eqref{size} and 
\be{infty-count} 
\| \sum_J \chi_{C2^s J}\|_{\infty} \lesssim 2^{Ns}
\end{equation}
hold. Let
$b_J$ be a collection of smooth functions satisfying \eqref{bj-prop}, and let 
$\psi_J$ be defined as above.
Then we have
\begin{equation}\label{next-reduction}
\| \sum_J \psi_J(b_J * S_{j+s} F_J) \|_p \lesssim 2^{-\eps s} 
(\sum_J |J| \|F_J\|_2^2)^{1/2}
\end{equation}
for all functions $F_J$ in $L^2(\H)$.
\end{proposition}

In this section we show how Proposition \ref{next} can be used to imply
Proposition \ref{first}.

Suppose that we are in the situation of Proposition \ref{first}.
We first observe a useful lemma which will also be needed much later in
this argument.

\begin{lemma}\label{bmo-mult} Let $B \subset B(0,C)$ be any 
Euclidean ball of side-length at least $2^{-\eps s}$, and define the functions 
$\psi_{J,B}$ by
$$ \psi_{J,B}(x) = \psi_B(2^{-j-s} \circ (x_J^{-1} x))$$
where $\psi_B$ is any bump function which is adapted to $B$.  Then we have
$$ | \{ \sum_J \psi_{J,B}(x) > s^3 2^{Ns} |B|\}| \lesssim 2^{-\eps s^2}.$$
\end{lemma}

\begin{proof} It suffices to show the two estimates
\be{l1-count} \| \sum_J \psi_{J,B}\|_1 \lesssim 2^{Ns} |B|
\end{equation}
and
\be{bmo-count} \| \sum_J \psi_{J,B}\|_{BMO} \lesssim s 2^{Ns} |B| 
\end{equation}
where BMO is defined with respect to the ball structure of the homogeneous group.  The desired distributional estimate follows from \eqref{l1-count} and \eqref{bmo-count} thanks to
the inequality
$$
| \{ |f| \geq \alpha \} \lesssim e^{-C\alpha/\|f\|_{BMO}} \frac{\|f\|_1}{\alpha}
$$
which follows immediately from the John-Nirenberg inequality and the Calder\'on-Zygmund
decomposition.

The estimate \eqref{l1-count} follows trivially from the triangle inequality
and \eqref{size}.  To show \eqref{bmo-count},
It suffices to show that
$$ \sum_J \osc_I \psi_{J,B} \lesssim 2^{Ns} |B|$$
for all balls $I$, where
$\osc_I f = \frac{1}{|I|} \int_I |f - f_I|$ and $f_I$ is the mean of $f$ on $I$.

Fix $I$, and suppose that $I$ has radius $2^i$. 
We divide into three cases, depending on the relative sizes of $I$, $J$ and
$2^s J$.  

We first consider the case where $I$ is larger than $2^s J$.  In this case
 $\osc_I \psi_{J,B}$ vanishes unless
$J$ is in $CI$, in which case the oscillation is $O(2^{Ns}|J| |B|/|I|)$.
Since the $J$ are disjoint and live in $CI$ the total contribution from these balls is acceptable.

Next, we consider the case where $I$ has size between $J$ and $2^s J$ inclusive.
For each scale $j$, there are at most $O(2^{Ns}|B|)$ balls $J$ of size $2^j$ which
give a non-zero contribution.  Since each ball contributes at most $O(1)$, we
are done.

Finally, we consider the case where $I$ has size smaller than $J$.  For
each scale $j$, there are at most $O(2^{Ns}|B|)$ balls $J$ which contribute.
But from the smoothness of $\psi_{J,B}$ we see that each ball gives
a contribution of $O(2^{-\eps(j+s-i)})$ for some $\eps > 0$.  Summing in $j$
we see that this contribution is also acceptable.
\end{proof}

By applying this lemma with a ball of size roughly 1 and a non-negative
cutoff, we obtain
\be{bigmult}
| \{ \sum_J \chi_{C 2^s J} \gtrsim s^3 2^{Ns}\}| \lesssim 2^{-\eps s^2}.
\end{equation}

To pass from this to \eqref{infty-count} we shall use a sieving argument
of C\'ordoba \cite{cordoba:sieve}.

For any ball $J \in \J$, define the \emph{height} $h(J)$ to be the number
$$ h(J) = \# \{ J' \in \J: 2J \subset 2J' \}.$$
We first deal with the contribution of those balls in $\J$ with height
at least $s^3 2^{Ns}$.  Clearly, the counting function
$\sum_{J \in \J} \chi_{C2^s J}$ is at least $s^3 2^{Ns}$ on these balls.
By the above lemma, the total measure of these balls is
$O(2^{-\eps s^2})$.  This implies that the contribution of these balls to
\eqref{first} is supported on a set of measure $O(2^{Ns} 2^{-\eps s^2})$,
which can safely be placed in the exceptional set $E$.

We now consider for each $a = 0, 1, \ldots, s^3-1$ the contribution
of those balls in $\J$ of height between $a 2^{Ns}$ and $(a+1) 2^{Ns}$.
If we denote this collection of balls by $\J_a$, then we claim that $\J_a$
obeys \eqref{infty-count}.  The estimate \eqref{first-reduction} would
then follow from $s^3$ applications of \eqref{next-reduction} and 
the triangle inequality.

It remains to verify \eqref{infty-count}.  Let $x$ be an arbitrary point
and let $\J^x$ be the set of all $J \in \J_a$ for which
$x \in 2^s J$.  We wish to show that $\# \J^x \lesssim 2^{Ns}$.

We may of course assume that $\J^x$ is non-empty.  Let $J_0$, $J_1$ be elements of
$\J^x$ with minimal and maximal radius $2^{j_0}$ and $2^{j_1}$ respectively.   
We observe that there are at most $O(2^{Ns})$ balls in $\J^x$ of 
radius comparable to $2^{j_0}$ since the balls are disjoint.  
Similarly there are
at most $O(2^{Ns})$ balls in $\J^x$ of radius comparable to $2^{j_1}$.  So it only
remains to show that there are at most $O(2^{Ns})$ balls in $\J^x$ of radius much larger than $2^{j_0}$ and much smaller than $2^{j_1}$.  But 
each such ball makes a contribution of $1$ to 
$h(j_1) - h(j_0)$, which is $O(2^{Ns})$ by assumption.
This completes the derivation of Proposition \ref{first} from Proposition
\ref{next}.

\section{Proof of Theorem \ref{weak-11} continued.  Iterated $TT^*$ methods.}

Let $\J$, $b_J$ satisfy the conditions of Proposition \ref{next}.  To finish
the proof of Theorem \ref{weak-11} it suffices
to show \eqref{next-reduction}.  By duality, it suffices to show
$$
(\sum_J |J|^{-1} \| S_{j+s}^*( \widetilde{b_J} * (\psi_J F) ) \|_2^2)^{1/2}
\lesssim 2^{-\eps s} \|F\|_{p'}$$
for all test functions $F$ on $\H$, where $\widetilde{b_J}(x) = b_J(x^{-1})$.
By the $TT^*$ method, it therefore suffices to show that
$$
\| \sum_J |J|^{-1} \psi_J(b_J * S_{j+s} S_{j+s}^* (\widetilde{b_J} * (\psi_J F)))\|_{p} \lesssim 2^{-\eps s} \|F\|_{p'}.$$
Since $S_{j+s} S_{j+s}^* = 2^{-Ns} |J|^{-1} S_0 S_0^*$, we may rewrite
this as
\be{new-targ}
\| TF \|_{p} \lesssim 2^{-\eps s} \|F\|_{p'},
\end{equation}
where 
$$ T = 2^{-Ns} \sum_J \psi_J T_J \psi_J$$
and $T_J$ is the self-adjoint operator
$$ T_J F = \frac{b_J}{|J|} * S_0 S_0^* (\frac{\widetilde{b_J}}{|J|} * F)).$$
Define the smooth functions $c_J$ supported on the ball $B(0,C)$ by
$$ c_J(v) = |J|^{-1} b_J(d_J(v))$$
where $d_J: B(0,C) \to CJ$ is the map
\be{dj-def}
 d_J(v) = x_J (2^j \circ v);
\end{equation}
from \eqref{bj-prop} we see that 
\be{cj-prop}
\|c_J\|_{L^1(B(0,C))} \lesssim 1, \quad \int_{B(0,C)} c_J = 0.
\end{equation}
Also, note that
$$S_0 S_0^* F(x) = \int \tilde\varphi(t) F(t \circ x)\ dt$$
where $\tilde\varphi$ is a bump function adapted to $\{t \sim 1\}$.  
We may rewrite $T_J$ as
\be{aj-def}
 T_J F(x) = 
\int \int \int c_J(v) \tilde\varphi(t) c_J(w) F(d_J(w) t \circ (d_J(v)^{-1} x))\ 
dw dt dv.
\end{equation}

We now define a slightly larger, and non-cancellative, version of $T_J$.
For each $J$, let $\psi^+_J$ be a slight enlargement of $\psi_J$
which is positive on the support of $\psi_J$.  Also,
apply Lemma \ref{integ} to find functions $c_J^1, \ldots, c_J^n$ supported
on $B(0,C)$ such that
\be{div}
 c_J = \sum_{i=1}^n \partial_{x_i} c_J^i
\end{equation}
and $\|c_J^i\|_1 \lesssim 1$.  If one then defines
$$ c_J^+ = |c_J| + \sum_{i=1}^n |c_J^i|,$$
then one sees that $c_J^+$ is a non-negative function on $B(0,C)$
with 
\be{cjp-prop}
\| c_J^+\|_1 \lesssim 1.
\end{equation}
Finally, we choose $\varphi^+$ to be any enlargement of $\tilde\varphi$ 
which is strictly positive on the support of $\tilde\varphi$, and obeys
the condition $\varphi^+(t) = t^{-N} \varphi^+(t^{-1})$.  
We then define the self-adjoint operator $T_J^+$ by
$$ T_J^+ F(x) = 
\int \int \int c^+_J(v) \varphi^+(t) c^+_J(w) F(d_J(w) t \circ (d_J(v)^{-1} x))\ dw 
dt dv,$$
and
$$ T^+ = 2^{-Ns} \sum_J \psi^+_J T^+_J \psi^+_J.$$
Clearly we have the pointwise bounds $T_J F(x) \leq T_J^+ F(x)$ and
$T F(x) \leq T^+ F(x)$ for all $J$ and non-negative $F$.

The operator
$$ F(\cdot) \mapsto F(d_J(w) t \circ (d_J(v)^{-1} \cdot))$$
is bounded on every $L^p$ uniformly in all variables.  From 
this, \eqref{aj-def}, \eqref{cjp-prop} and
Minkowski's inequality we 
therefore have $\| T_J^+ F \|_p \lesssim \|F\|_p$ uniformly in $p$ and $J$.

We now show
\be{ap-bound}
\| T^+ F \|_{p} \lesssim \|F\|_q
\end{equation}
for all $1 \leq p \leq q \leq \infty$; note that this would imply \eqref{new-targ} were it not for the $2^{-2\eps s}$ factor.  By interpolation and duality it
suffices to verify this for $q = \infty$.  Since $T_J^+$ is positivity
preserving it suffices to verify this when $F$ is identically 1, i.e. we need
to show that
$$
\| 2^{-Ns} \sum_J \psi_J T_J^+ \psi_J \|_{p} \lesssim 1.
$$
By \eqref{infty-count} and H\"older's inequality we have the pointwise estimate
$$ 2^{-Ns} \sum_J \psi_J T_J^+ \psi_J
\lesssim 2^{-Ns/p} (\sum_J (T_J^+ \psi_J)^p)^{1/p}.$$
Thus to show \eqref{ap-bound} it suffices to show
$$ 2^{-Ns/p} (\sum_J \|T_J^+ \psi_J\|_p^p)^{1/p} \lesssim 1.$$
But this follows from the $L^p$ boundedness of $T_J^+$ and
\eqref{size}.

To finish the proof of Theorem \ref{weak-11} we must obtain the
gain of $2^{-2\eps s}$ for \eqref{new-targ}.
To obtain this gain, we will iterate $T$ $m$ times as in the proof of 
Theorem \ref{L2}, until the kernel is smooth enough to profitably
interact with the derivatives in \eqref{div}.  As long as there is
enough isotropic smoothing, we hope to bound (in an appropriate
sense) $T^m$ by $2^{-\eps s} (T^+)^m$.  As before, there will be an 
exceptional portion of $T^m$, but we hope to also show this part
is also very small.

Naively, one expects this to work with $m=n$.  However, one runs into
two difficulties with this choice.  Firstly, if one composes
the operator $T_I$ with $T_J$, where $I$ is much larger than $J$,
the cutoff functions $\psI$ corresponding to $J$ can
seriously truncate the smoothing effect from $T_I$.
  Secondly, when $\H$ is a general
homogeneous group, the smoothing effects of the $T_J$
will tend to be along almost parallel directions, rather than being
isotropically dispersed.  Although each of these obstacles is individually
tractable, the combined effect of these two obstacles may restrict
the smoothing effects discussed above to very short, very
parallel arcs, which will not give much isotropic 
regularity\footnote{
These obstructions also show that in the non-Euclidean case
the smoothing effect is global rather than local, and one cannot  
exhibit this effect by naive microlocal 
methods (as used in the standard proof of the averaging lemma); this 
is in contrast with the Fourier transform analysis of
\cite{seeger:rough} in the Euclidean case.}.
To avoid
this problem we
shall iterate $T$ by considerably more than $n$ times
to ensure the existence of at least $n$ untruncated arcs.  In fact
we shall iterate $m = 2^{2n - 3}$ times.

We now turn to the details.  Let $m$ be a large number to be chosen later.  To show \eqref{new-targ}
it suffices to show that
\be{am-bound}
 \| T^m F \|_p \lesssim 2^{-\eps s} \| F\|_\pp
\end{equation}
for some $\eps > 0$.  To see this, observe from the $TT^*$ method and
the self-adjointness of $T$ that
\eqref{am-bound} implies
$$ 
\| T^{m/2} F \|_p \lesssim 2^{-\eps s} \|F\|_2$$
with a slightly worse value of $\eps$.  On the other hand, from 
many applications of \eqref{ap-bound} we have
$$ 
\| T^{m/2} F \|_p \lesssim \|F\|_q$$
for all $q \geq p$.  By interpolation we thus obtain
$$ 
\| T^{m/2} F \|_p \lesssim 2^{-\eps s} \|F\|_\pp$$
for an even worse value of $\eps$.  By iterating this argument $2n-3$ times we thus obtain \eqref{new-targ} (for a very small
value of $\eps$).

It remains to show \eqref{am-bound}.  Since $A^m$ is bounded on $L^2$
by \eqref{ap-bound}, it suffices by interpolation
to prove this for $p=1$.  By expanding $T^m$, we thus reduce
to showing that
$$
 \|2^{-Nms} \sum_{J_1 \ldots J_m \in \J} (\prod_{i=1}^m \psi_{J_i} T_{J_i} \psi_{J_i}) F \|_1 \lesssim 2^{-\eps s} \|F\|_\infty.
$$
We use $2^{j_i}$ to denote the radius of $J_i$.

The balls $J_i$ may be radically different sizes, and need not be arranged
in any sort of monotone order.  Nevertheless, we can still extract a 
subsequence of $n$ balls whose sizes do increase monotonically, and have
no smaller balls between elements of the sequence.  More precisely, we have

\begin{definition} Let $k = (k_1, \ldots, k_n)$ be a strictly increasing
$n$-tuple of integers in $\{1, \ldots, m\}$.  We say that an $m$-tuple $J = (J_1, \ldots, J_m)$ of balls is \emph{ascending} with respect to $k$ if
$$ j_{k_q} \leq j_l \hbox{ for all } k_q \leq l \leq k_n,$$
and write this as $J \nearrow k$.  Similarly, we say that $J$ is \emph{descending} with respect to $k$ if
$$ j_{k_q} \leq j_l \hbox{ for all } k_1 \leq l \leq k_q,$$
and write this as $J \searrow k$.
\end{definition}

\begin{lemma}
If $m \geq 2^{2n-3}$ and $J \in \J^m$, then there exists a sequence $k$
such that either $J \nearrow k$ or $J \searrow k$.
\end{lemma}

\begin{proof}  We first construct an auxilliary sequence $l_1, \ldots,
l_{2n-2}$ of integers and a sequence $S_1, \ldots S_{2n-2}$ of intervals
of integers by the following iterative procedure.  

Let $S_1$ be the interval $\{1, \ldots, m\}$.  
For each $p = 1, \ldots, 2n-2$ in turn, we
choose $l_p \in S_p$ so that $J_{l_q}$ has minimal radius among all
the balls $\{J_l: l \in S_p\}$.  Removing the element $l_p$ from $S_p$
divides the remainder into two intervals $\{l \in S_p: l < l_p\}$ and $\{l \in
S_p: l > l_p\}$; we choose $S_{p+1}$ to be the larger of the two
intervals.  We then increment $p$ and iterate the above construction.

One can easily show inductively that $|S_p| \geq 2^{2n-2-p}$ for all
$s$, so that all the $l_p$ are well defined.  Furthermore, one has
$j_{l_q} \leq j_l$ for all $l$ between $l_q$ and $l_{2n-2}$.
One of the sets $\{p: l_p \leq l_{2n-2}\}$, $\{p: l_p \geq l_{2n-2}\}$
has a cardinality of at least $n$.  If the former set is larger, we
choose $k$ to be the first $n$ elements of this set and observe that
$J \nearrow k$.  Otherwise we choose $k$ to be the first $n$ elements
of the latter set and observe that $J \searrow k$.
\end{proof}

Temporarily set $m = 2^{2n-3}$.
Order the sequences $k$ lexicographically, so in particular we
have $k < k'$ whenever $k_1 < k'_1$.  For all $J \in J^m$, let
$k_{max}(J)$ be the largest sequence with respect to this ordering so that
either $J \nearrow k$ or $J \searrow k$.  From the above lemma we 
see that $k_{max}(J)$ is well-defined.  Since the number of sequences is finite,
it suffices to show that
$$\|2^{-Nms} \sum_{J_1, \ldots, J_m \in \J: k_{max}(J_1, \ldots, J_m) = k} (\prod_{i=1}^m \psi_{J_i} T_{J_i} \psi_{J_i}) F \|_1 \lesssim 2^{-\eps s} \|F\|_\infty$$
for each $k$.

Fix $k$.  The purpose of the following (somewhat technical) discussion is to enable us to reduce to the case when $k_1 = 1$ and $k_n = m$.

We observe from the lexicographical ordering that the property that $k_{max}(J_1, \ldots, J_m) = k$ is independent of the choices of $J_i$ for $1 \leq i < k_1$.  We thus abuse notation and
write
$$ k_{max}(J_{k_1}, \ldots, J_m) = k \hbox{ instead of } k_{max}(J_1, \ldots, J_m) = k.$$
The desired estimate can then be factored as
$$\|2^{-N(m-k_1+1)s} \sum_{J_{k_1}, \ldots, J_m \in \J: k_{max}(J_{k_1}, \ldots, J_m) = k} (\prod_{i=k_1}^m \psi_{J_i} T_{J_i} \psi_{J_i}) T^{k_1-1} F \|_1 \lesssim 2^{-\eps s} \|F\|_\infty.$$
By \eqref{ap-bound}, $T$ is bounded on $L^\infty$, and it suffices
to show that
$$2^{-N(m-k_1+1)s} \|\sum_{J_{k_1}, \ldots, J_m \in \J: k_{max}(J_{k_1}, \ldots, J_m) = k} (\prod_{i=k_1}^m \psi_{J_i} T_{J_i} \psi_{J_i}) F \|_1 \lesssim 2^{-\eps s} \|F\|_\infty.$$
The left-hand side is majorized by
\be{tmp}
2^{-N(m-k_1+1)s} 
\| \sum_{J_{k_1}, \ldots, J_m \in \J: J \nearrow k \hbox{ or } J \searrow k } 
|(\prod_{i=k_1}^m \psi_{J_i} T_{J_i} \psi_{J_i}) F|\|_1
\end{equation}
where $J = (J_1, \ldots, J_m)$, and the choices of $J_1, \ldots, J_{k_1-1}$ are irrelevant.  Since the properties $J \nearrow k$, $J \searrow k$ do not
depend on $J_{k_n+1}, \ldots, J_m$, we may estimate \eqref{tmp} crudely by
$$
2^{-N(k_n-k_1+1)s} 
\| (T^+)^{m-k_n} \sum_{J_{k_1}, \ldots, J_{k_n} \in \J: J \nearrow k \hbox{ or } J \searrow k } |(\prod_{i=k_1}^{k_n} \psi_{J_i} T_{J_i} \psi_{J_i}) F|\|_1
$$
By \eqref{ap-bound} we may discard the $(T^+)^{m-k_n}$ operator.
By a re-labelling of $J$ and $k$, and reducing $m$ to $k_n-k_1+1$, it thus suffices to show that
$$ 2^{-Nms} 
\| \sum_{J_{1}, \ldots, J_{m} \in \J: J \nearrow k \hbox{ or } J \searrow k } |(\prod_{i=1}^{m} \psi_{J_i} T_{J_i} \psi_{J_i}) F|\|_1
\lesssim 2^{-\eps s} \|F\|_\infty
$$ 
for all $m \leq 2^{2n-3}$ and all $k$ such that $k_1 = 1$, $k_n = m$.

Fix $m$, $k$.  By duality it suffices to show that
\be{bilinear}
2^{-Nms} 
\sum_{J_{1}, \ldots, J_{m} \in \J: J \nearrow k \hbox{ or } J \searrow k } |\langle (\prod_{i=1}^{m} \psi_{J_i} T_{J_i} \psi_{J_i}) F_J, G_J \rangle|
\lesssim 2^{-\eps s} 
\end{equation} 
for all functions $F_J$, $G_J$ in the unit ball of $L^\infty$.  It suffices
to consider the contribution of $J \nearrow k$, since the other
contribution then follows by self-adjointness.

For each $J \nearrow k$, we expand the inner product in \eqref{bilinear} as
\be{expand}
\int \int \int \int G_J(x_0) F_J(x_m) (\prod_{i=1}^m
\psi_{J_i}(x_{i-1}) c_{J_i}(v_i) \tilde\varphi(t_i) c_{J_i}(w_i)
\psi_{J_i}(x_i))\ dx_0 dw dt dv
\end{equation}
where $v = (v_1, \ldots, v_n)$, $w = (w_1, \ldots, w_n)$ range over $B(0,C)^n$, 
$t = (t_1, \ldots, t_n)$ ranges
over $[C^{-1},C]^n$, $dw = \prod_{i=1}^n dw_i$, $dv = \prod_{i=1}^n dv_i$,
$x_0$ ranges over $\H$, and $x_1, \ldots, x_m$ are defined recursively
by
\be{xi-def}
 x_i = F(d_{J_i}(w_i) t_i \circ (d_{J_i}(v_i)^{-1} x_{i-1})) \hbox{ for } i=1, \ldots, m.
\end{equation}
Note that each $x_i$ is a function of $x_0$ and $J_l, v_l, t_l, w_l$ for
all $l = 1, \ldots, i$.  We call the variables $x_0$ and
$v_l, t_l, w_l$ for $l=1, \ldots, m$ \emph{integration variables}.

There are many variables of integration here, but the only ones that
we shall actively use are the dilation parameters $t_{k_1}, \ldots, t_{k_n}$
and the translation parameter $v_1$.  Accordingly, we define new
variables $\tau = (\tau_1, \ldots, \tau_n), y = (y_1, \ldots, y_n) \in \R^n$ by $\tau_q = t_{k_q}$
and $y = v_1$.

Each $\tau_q$ integration is smoothing in one direction.
The combined smoothing effect of all the $\tau$ variables shall
be beneficial provided that the Jacobian
$$\det D^L_\tau(x_m)$$
is sufficiently large.

As will become clear later, the natural size for $\det D^L_\tau$ is $2^{M_n}$, where
the quantities $M_0, \ldots, M_n$ are defined by
\be{M-def}
M_q = \sum_{i=1}^q \alpha_i (j_{k_i}+s).
\end{equation}

Accordingly, we shall decompose the $J \nearrow k$ portion of \eqref{bilinear} into
\be{non-cancel-l1}
\bs
\sum_{J \in \J^m: J \nearrow k} 
|&\int\int\int\int G_J(x_0) F_J(x_m) (\prod_{i=1}^m
\psi_{J_i}(x_{i-1}) c_{J_i}(v_i) \tilde\varphi(t_i) c_{J_i}(w_i)
\psi_{J_i}(x_i))\\
& \eta(2^{\delta s} 2^{-M_n} \det D^L_\tau(x_m)) \ dx_0
dw dt dv|
\lesssim 2^{-\eps s} 2^{Nms}
\end{split}
\end{equation}
and
\be{cancel-l1}
\bs
\sum_{J \in \J^m: J \nearrow k} 
|&\int\int\int\int G_J(x_0) F_J(x_m)
(\prod_{i=1}^m
\psi_{J_i}(x_{i-1}) c_{J_i}(v_i) \tilde\varphi(t_i) c_{J_i}(w_i)
\psi_{J_i}(x_i)) \\
&(1-\eta(2^{\delta s} 2^{-M_n} \det D^L_\tau(x_m))) \ dx_0
dw dt dv|
\lesssim 2^{-\eps s} 2^{Nms}.
\end{split}
\end{equation}
Here $\delta > 0$ is a small number to be chosen later,
and $\eta$ is a bump function
which equals $1$ near 1.

We shall prove these two estimates in later sections.  But first we must
introduce some preliminaries to treat the Jacobian $\det D^L_\tau(x_m)$,
which is the wedge product of $n$ vectors of vastly different sizes.

\section{The exterior algebra and non-isotropic scaling}

It shall be necessary to define some artificial structures on
the exterior algebra $\Lambda$ of $\R^n$.
Define a quasi-order $\precsim$ on $\Lambda$ by 
$$ \sum_P a_P e_P \precsim \sum_P b_P e_P \quad \iff \quad
|a_P| \lesssim b_P \hbox{ for all } P,$$
and write $w \approx w'$ if $w \precsim w'$ and $w' \precsim w$.
We define an absolute value by
$$ \|\sum_P a_P e_P\| = \sum_P |a_P| e_P.$$
Of course, $\|a\| = |a|$ for scalars $a$.
We let $\one$ denote the vector $\one = (1, \ldots, 1)$.

We define a non-cancellative analogue 
$\diamond$ of the wedge product by 
$$ (\sum_P a_P e_P) \diamond (\sum_Q b_Q e_Q) = \sum_P \sum_Q a_P b_Q \|e_P \wedge e_Q\|.$$
Note that $\diamond$ is bilinear and associative.  This operation
dominates the wedge product in the following sense:
if $\omega_1 \precsim a_1$ and $\omega_2 \precsim a_2$ then
\be{wedge-order} \|\omega_1 \wedge \omega_2\| \precsim a_1 \diamond a_2, \quad
|\omega_1 \cdot \omega_2| \lesssim a_1 \cdot a_2.
\end{equation}
Finally, we observe that if
$1 \leq r \leq n$ and 
$i_1, \ldots, i_r$ is any non-decreasing sequence of integers, then
\be{diamond}
(2^{i_1} \circ \one) \diamond \ldots \diamond (2^{i_r} \circ \one)
\approx \sum_{1 \leq p_1 < \ldots < p_r \leq n}
2^{\alpha_{p_1} i_1 + \ldots + \alpha_{p_r} i_r} e_{p_1} \wedge \ldots \wedge
e_{p_r}.
\end{equation}
In particular, from \eqref{M-def} and the hypothesis $J \nearrow k$ we have
\be{m-size} (2^{j_{k_1}+s} \circ \one) \diamond \ldots \diamond (2^{j_{k_n}+s} \circ \one)
\approx 2^{M_n} e_1 \wedge \ldots \wedge e_n.
\end{equation}

If $F(x)$ is a form-valued function of $x$ and $C > 0$, define
$$ \| (1 + C \nabla_x)F(x) \| = \| F(x) \| + C \sum_{i=1}^n \| \partial_{x_i} F(x)\|.$$
More generally, we define
$$ \| (A + B \nabla_x + C \nabla_y) F(x,y) \| = A \|F(x,y)\|
+ B \sum_{i=1}^n \|\partial_{x_i} F(x,y)\| + C \sum_{q=1}^n \| \partial_{y_q} F(x,y) \|.$$
From the product rule and \eqref{wedge-order} we observe that
\be{product}
\| (1 + C\nabla_x) (F \cdot G) \| \lesssim \| (1 + C \nabla_x) F \| \cdot
\| (1 + C\nabla_x) G\|
\end{equation}
and
\be{wedgehog}
\| (1 + C\nabla_x) (F \wedge G) \| \precsim \| (1 + C \nabla_x) F \| \diamond
\| (1 + C\nabla_x) G\|.
\end{equation}

We record the following estimates on the size and derivatives of $x_l$, and $\det D^L_T(x_m)$.

\begin{lemma}\label{derivs}  If $t_l \sim 1$ and $x_l \in (2^s J_l)_\Delta$
for all $1 \leq l \leq m$, then
\begin{align}
\|(1 + \nabla_\tau) \partial^L_{y_i} x_l\| 
&\precsim 2^{j_{1} + s - cs} \circ \one 
\label{xl-yi}\\
\|(1 + 2^{cs} \nabla_y + \nabla_\tau) \partial^L_{\tau_q} x_m\| &\precsim 2^{j_{k_q} + s} \circ \one
\label{xl-tq-v}\\
\|(1 + 2^{cs} \nabla_y + \nabla_\tau) \det D^L_\tau(x_m)\| &\lesssim 2^{M_n}\label{jac-y}\\
\|(1 + 2^{cs} \nabla_y + \nabla_\tau) \psi_{J_l}(x_{l'})\| &\lesssim
\psi^+_{J_l}(x_{l'})\label{psi}
\end{align}
for all $i, q,q' = 1, \ldots n$, $1 \leq l \leq m$, $0 \leq l' \leq l$,
where $c>0$ is a constant independent of $\eps$, $\delta$.
\end{lemma}

\begin{proof}  
From \eqref{xi-def}, \eqref{dj-def},
\eqref{group-diff}, and \eqref{dil-diff} we have
\be{yi-form} \partial^L_{y_i} x_l = (t_l \ldots t_0)
\circ C[x_{J_1}^{-1} x_0] \partial^L_{y_i} (2^{j_1} \circ y^{-1}).
\end{equation}
By \eqref{c-twine} and \eqref{dil-diff}, this becomes
$$ \partial^L_{y_i} x_l = (t_l \ldots t_0 2^{j_1 + s})
\circ C[2^{-j_1 - s} \circ x_{J_1}^{-1} x_0] (2^{-s} \circ \partial^L_{y_i} y^{-1}).$$
Since $x_0 \in (2^s J_1)_\Delta$, $2^{-j_1 - s} \circ x_{J_1}^{-1} x_0$ is bounded.
From this, \eqref{c-bound} and the observation that
$|\partial^L_{y_i} y^{-1}| \lesssim 1$, we thus have
$$ |C[2^{-j_1 - s} \circ x_{J_1}^{-1} x_0] (2^{-s} \circ \partial^L_{y_i} y^{-1})|
\lesssim 2^{-cs}.$$
Inserting this into the previous estimate we thus obtain
$$
 \partial^L_{y_i} x_l \precsim 2^{j_1 + s - cs} \circ \one
$$
which is the first part of \eqref{xl-yi}.  The $\nabla_\tau$ portion of
\eqref{xl-yi} then follows from \eqref{yi-form} and \eqref{trivx-bound}.

We now turn to \eqref{xl-tq-v}.  From \eqref{xi-def}, \eqref{group-diff},
and \eqref{dil-diff}, we have
\be{vi-formula}
\partial^L_{\tau_q} x_m
= t_{k_q}^{-1} ((t_m \ldots t_{k_q} 2^{j_{k_q} + s}) \circ X(u_q))
\end{equation}
where $u_q$ is the quantity
\be{uq-def}
u_q = 2^{-j_{k_q} - s} \circ (d_{J_{k_q}}(v_{k_q})^{-1} x_{k_q-1}))
= (2^{-j_{k_q}-s} \tau_q^{-1}) \circ (d_{J_{k_q}}(w_{k_q})^{-1} x_{k_q}).
\end{equation}
Since $x_{k_q-1} \in (2^s J_{k_q})_\Delta$, $u_q$ is bounded,
and so the first part of \eqref{xl-tq-v} obtains.  To show the $2^{cs} \nabla_y$
portion of \eqref{xl-tq-v}, it suffices from \eqref{vi-formula} and the chain rule
to show that
$\|\partial^L_{y_i} u_q \| \lesssim 2^{-cs}$.
But by \eqref{uq-def}, \eqref{group-diff}, and \eqref{dil-diff} we have
$$ \partial^L_{y_i} u_q = (\tau_q^{-1} 2^{-j_{k_q}-s}) \circ \partial^L_{y_i} x_{k_q},$$
and the claim follows from \eqref{xl-yi} and the inequality
$j_{k_1} \leq j_{k_q}$ arising from the hypothesis $J \nearrow k$.

We now show the $\nabla_\tau$ portion of \eqref{xl-tq-v}.  We consider the
$\partial_{\tau_{q'}}$ derivatives for $q' \geq q$ and $q' < q$ separately.
If $q' \geq q$, then we see from \eqref{vi-formula} and \eqref{trivx-bound}
that
$$ \rho(\partial_{\tau_{q'}} \partial^L_{\tau_q} x_m)
\lesssim \rho(\partial^L_{\tau_q} x_m),$$
so the claim follows from the first part of \eqref{xl-tq-v}.  If $q' < q$,
then from \eqref{vi-formula} we have
$$
\partial_{\tau_{q'}} \partial^L_{\tau_q} x_m = t_{k_q}^{-1} (t_m \ldots t_{k_q} 2^{j_{k_q} + s}) \circ \partial_{\tau_{q'}} X(u_q).$$
Since $X$ is polynomial and $u_q$ is bounded, it thus suffices by \eqref{comparable}
to show that
$$ |\partial^L_{\tau_{q'}} u_q| \lesssim 1.$$
But from \eqref{uq-def}, \eqref{group-diff}, and \eqref{dil-diff}, we have
$$ \partial^L_{\tau_{q'}} u_q = 2^{-j_{k_q} - s} \circ \partial^L_{\tau_{q'}} x_{k_q-1},$$
and the claim follows from the first part of \eqref{xl-tq-v}.

We now turn to \eqref{jac-y}.  It suffices to show that
$$
|(1 + 2^{cs} \nabla_y + \nabla_\tau) (\partial^L_{\tau_1} x_m
\wedge \ldots \wedge \partial^L_{\tau_n} x_m)|
\precsim 2^M e_1 \wedge \ldots \wedge e_n.$$
From \eqref{wedgehog} and \eqref{wedge-order} we have
\bas
|(1 + 2^{cs} \nabla_y + \nabla_\tau) &(\partial^L_{\tau_1} x_m
\wedge \ldots \wedge \partial^L_{\tau_n} x_m)|\\
&\precsim
|(1 + 2^{cs} \nabla_y + \nabla_\tau) \partial^L_{\tau_1} x_m|
\diamond \ldots \diamond
|(1 + 2^{cs} \nabla_y + \nabla_\tau) \partial^L_{\tau_n} x_m|.
\end{align*}
By \eqref{xl-tq-v} this is majorized by
$$ (2^{j_{k_1}+s} \circ \one) \diamond \ldots \diamond (2^{j_{k_l}+s} \circ \one).$$
The claim then follows from \eqref{m-size}.

Finally, we show \eqref{psi}.
We can rewrite the desired estimate as
$$
|(1 + 2^{cs} \nabla_y + \nabla_\tau) \psi(2^{-j_l-s} \circ (x_{J_l}^{-1} x_{l'}))| 
\lesssim \psi^+(2^{-j_l-s} \circ (x_{J_l}^{-1} x_{l'})).$$
From the support assumptions on $\psi$ and $\psi^+$ we have $|(1 + \nabla)\psi|
\lesssim \psi^+$.  Thus by the chain rule and \eqref{comparable}, it suffices to show that
$$ |2^{cs} \partial^L_{y_i} (2^{-j_l-s} \circ (x_{J_l}^{-1} x_{l'})| \lesssim 1$$
and
$$ |\partial^L_{\tau_q} (2^{-j_l-s} \circ (x_{J_l}^{-1} x_{l'})| \lesssim 1$$
for all $i,q = 1, \ldots, n$.  We may of course assume that $1 \leq l'$ and
$k_q \leq l'$ since the claims are trivial otherwise.
But these estimates follow from \eqref{group-diff},
\eqref{dil-diff}, \eqref{xl-yi}, and \eqref{xl-tq-v}, noting that
$j_1, j_{k_q} \leq j_{l}$ from the hypothesis $J \nearrow k$.
\end{proof}

\section{Proof of Theorem \ref{weak-11} continued.  The degenerate portion
of the integral.}

We now prove \eqref{non-cancel-l1}.  For this estimate we do not
exploit any cancellation, and crudely majorize the left-hand side
as
\be{non-deg}
\sum_{J \in \J^m: J \nearrow k}
\int_{|\det D^L_\tau x_m| \lesssim 2^{-\delta s} 2^{M_n}} \prod_{i=1}^m
\psi^+_{J_i}(x_{i-1}) 
c^+_{J_i}(v_i) \varphi^+(t_i) c^+_{J_i}(w_i)
\psi^+_{J_i}(x_i)\ dx_0
dw dt dv
\end{equation}
We discard
the $\psi^+_{J_i}(x_i)$ multiplier.  We may freeze the
$t_i$, $v_i$, $w_i$ variables using \eqref{cjp-prop}
and reduce ourselves to showing
\be{nondeg-targ}
\sum_{J \in \J^m: J \nearrow k}
\int_{|\partial^L_{\tau_1} x_m \wedge \ldots \wedge \partial^L_{\tau_n} x_m| \lesssim 2^{-\delta s} 2^{M_n}} 
\prod_{i=1}^m
\psi^+_{J_i}(x_{i-1})
\ dx_0
\lesssim 2^{-\eps s} 2^{Nms}
\end{equation}
uniformly over all choices of $t_i \sim 1$, $v_i \in B(0,C)$,
$w_i \in B(0,C)$, where $v_i$ and $w_i$ are allowed to depend on $J_i$.

Fix $t$, $w$, $v$.  To show \eqref{nondeg-targ},
we first exclude an exceptional set of $x$'s.

\begin{definition}  For each $x$ in $\H$, define the set $S(x)$
by
$$ S(x) = \{ 2^{-j-s} \circ (x_J^{-1} x): J \in \J, x \in (2^s J)_\Delta \}.$$
A point $x$ is said to be \emph{good} if one has
\be{equi} \# S(x) \cap B \lesssim s^3 2^{Ns} |B|
\end{equation}
for all balls $B$ of radius $2^{-\eps s}$.
\end{definition}

From \eqref{infty-count} we see that $S(x)$ is supported in
the unit annulus $A_0$ and $\# S(x) \lesssim 2^{Ns}$ for all $x$.  The
property \eqref{equi} can thus be thought of as a statement about the
uniform distribution of $S(x)$.

Let $E$ denote the set of all points in $x$ which are not good.
Fortunately, $E$ is very small:

\begin{lemma}  If $\eps$ is sufficiently small, we have
\be{e-size}
|E| \lesssim 2^{-\eps s^2}.
\end{equation}
\end{lemma}

\begin{proof}  Since there are at most $O(2^{C s})$ 
finitely overlapping $B$ which need to be considered for \eqref{equi}, it suffices
to show that
$$ | \{ x \in \H: \# S(x) \cap B \gtrsim s^3 2^{Ns} |B| \}| \lesssim 2^{-\eps s^2}$$
for each ball $B$.  But this follows from Lemma \ref{bmo-mult} after
some re-arranging.
\end{proof}

For each $i = 1,\ldots,m$,
the contribution to \eqref{non-deg} of the case when $x_i \in E$
is bounded by
$$\sum_{J \in \J^m}
\int
\prod_{l=1}^n \psi^+_{J_l}(x_{l-1})
\chi_E(x_l)
\ dx_0.$$
By \eqref{infty-count} this is bounded by
$$ 2^{Nms} \int \chi_E(x_l)\ dx_0.$$
Thus the contribution to \eqref{non-deg} is definitely acceptable
by \eqref{e-size} and the observation that
$x_0 \mapsto x_l$ is a diffeomorphism with Jacobian
\be{diffeo}
\det D_{x_0}(x_{l}) = (t_1 \ldots t_l)^N \sim 1
\end{equation}

Thus it remains only to show that
$$
\sum_{J \in \J^m: J \nearrow k}
\int_{|\partial^L_{\tau_1} x_m \wedge \ldots \wedge \partial^L_{\tau_n} x_m| \lesssim 2^{-\delta s} 2^{M_n}} 
\prod_{i=1}^m
\psi^+_{J_i}(x_{i-1})
\chi_{E^c}(x_{i-1})
\ dx_0
\lesssim 2^{-\eps s} 2^{Nms}.
$$

For each $q=0,\ldots,n$, define $P_q$ to be the property that
$$ 2^{-q\delta s/n} 2^{M_q} \lesssim | \partial^L_{\tau_1} x_m
\wedge \ldots \wedge \partial^L_{\tau_q} x_m \cdot
e_1 \wedge \ldots \wedge e_q|,$$
where $M_q$ was defined in \eqref{M-def}. The desired estimate can thus be rewritten as
$$
\sum_{J \in \J^m: J \nearrow k}
\int_{P_n \hbox{ fails}} 
\prod_{i=1}^m
\psi^+_{J_i}(x_{i-1})
\chi_{E^c}(x_{i-1})
\ dx_0
\lesssim 2^{-\eps s} 2^{Nms}.
$$

Since $P_0$ is vacuously true, it thus suffices to show
\be{non-deg-targ}
\sum_{J \in \J^m: J \nearrow k}
\int_{P_{q-1} \hbox{ holds}, P_q \hbox{ fails}}
\prod_{i=1}^m
\psi^+_{J_i}(x_{i-1})
\chi_{E^c}(x_{i-1})
\ dx_0
\lesssim 2^{-\eps s} 2^{Nms}
\end{equation}
for all $q=1,\ldots, n$ (cf. \eqref{ei-decomp}).  

Fix $1 \leq q \leq n$.  We now make the key observation

\begin{proposition}  If we fix $x_0$ and all the $J_i$ except for
$J_{k_q}$, then we have
\be{jkq-card}
 \# \{ J_{k_q} : J \in \J^m_k, P_{q-1} \hbox{ holds}, P_q \hbox{ fails} \}
\lesssim 2^{-\eps s} 2^{Ns}
\end{equation}
provided that $x_{k_q-1}$ is good.
\end{proposition}

\begin{proof} 
Suppose $J_{k_q}$ is in the set in \eqref{jkq-card}.  Since $P_q$ fails, we have
$$ 
| (\partial^L_{\tau_1} x_m
\wedge \ldots \wedge \partial^L_{\tau_q} x_m) \cdot 
(e_1 \wedge \ldots \wedge e_q) | \lesssim 2^{-q\delta s/n} 2^{M_q}.$$
We rewrite this as
\be{vec}
 | (2^{-j_{k_q}-s} \circ \partial^L_{\tau_q} x_m) \cdot a | \lesssim 2^{-q\delta s/n}
\end{equation}
where the vector $a = a_1 e_1 + \ldots + a_q e_q$ is defined by
$$
a_l = 2^{\alpha_l (j_{k_q}+s)} 
2^{-M_q} (\partial^L_{\tau_1} x_m
\wedge \ldots \wedge \partial^L_{\tau_{q-1}} x_m \wedge e_l) \cdot 
(e_1 \wedge \ldots \wedge e_q).
$$
Since $M_q = M_{q-1} + \alpha_q (j_{k_q}+s)$, we may rewrite this as
\be{al-def} 
a_l = \pm 2^{-(\alpha_q - \alpha_l)(j_{k_q}+s)}
2^{-M_{q-1}} (\partial^L_{\tau_1} x_m
\wedge \ldots \wedge \partial^L_{\tau_{q-1}} x_m) \cdot
(e_1 \wedge \ldots \widehat{e_l} \ldots \wedge e_q)
\end{equation}
where the $\widehat{e_l}$ denotes that the $e_l$ term is missing from the wedge
product.
Since $P_{q-1}$ holds, we thus see that 
\be{large} |a_q| \gtrsim 2^{-(q-1)\delta s/n}.
\end{equation}

Also, from \eqref{xl-tq-v} we see that
$$ |a_l| \lesssim 2^{-(\alpha_q - \alpha_l)(j_{k_q}+s)}
2^{-M_{q-1}} ( 
(2^{j_{k_1}+s} \circ \one) \diamond \ldots \diamond (2^{j_{k_{q-1}}+s} \circ \one)
) \cdot (e_1 \wedge \ldots \widehat{e_l} \ldots \wedge e_q).$$
By \eqref{diamond}, we thus have
$$ |a_l| \lesssim 2^{-(\alpha_q - \alpha_l)(j_{k_q}+s)}
2^{-M_{q-1}} 
(\prod_{l'=1}^{l-1} 2^{\alpha_{l'}(j_{k_{l'}}+s)}) (\prod_{l'=l}^{q-1}
2^{\alpha_{l'+1}(j_{k_{l'}}+s)}).$$
By \eqref{M-def}, this simplifies to
$$ |a_l| \lesssim 2^{-(\alpha_q - \alpha_l)(j_{k_q}+s)}
\prod_{l'=l}^{q-1}
2^{(\alpha_{l'+1}-\alpha_{l'})(j_{k_{l'}}+s)}).$$
Since $J \nearrow k$, we have $j_{k_{l'}} + s \leq j_{k_{q-1}} + s$.
Applying this inequality, we obtain a telescoping product which simplifies to
\be{al-est}
 a_l \lesssim 2^{-(\alpha_q - \alpha_l)(j_{k_q} - j_{k_{q-1}})}
\end{equation}

From \eqref{al-def} we see that $a_l$ is independent of $j_{k_q}$
if $\alpha_l = \alpha_q$.
If $\alpha_l < \alpha_q$, then $a_l$ can vary with $j_{k_q}$.  However, 
from \eqref{al-est} we see that $a_l = O(2^{-Cs})$ unless
$j_{k_q} = j_{k_{q-1}} + O(s)$.

In both cases we thus conclude that, up to an error of $2^{-Cs}$, 
the quantities $a_l$ can
each take at most $O(s)$ values.  From this, \eqref{large}, and \eqref{vec},
we see that
$2^{-j_{k_q}-s} \circ \partial^L_{\tau_q} x_m$ lies in a union of
$O(s^C)$ $O(2^{-\delta s/n})$-neighbourhoods of hyperplanes.

From \eqref{vi-formula} and the fact that the frozen quantities $t_i$
are comparable to 1, we thus see that $X(u_q)$
also lives in a union of $O(s^C)$ $O(2^{-\delta s/n})$-neighbourhoods
of hyperplanes. 
From Lemma \ref{x-invert} and the boundedness of $u_q$,
we thus see that $u_q$ lives in a union
of $O(s^C)$ $O(2^{-\delta s/n})$-neighbourhoods of compact hypersurfaces.
From \eqref{uq-def} we have
$$ u_q = (2^{-s} \circ v_{k_q}^{-1}) 2^{-j_{k_q}-s} \circ (x_{J_{k_q}}^{-1} x_{k_q-1}),$$
and so $2^{-j_{k_q}-s} \circ (x_{J_{k_q}}^{-1} x_{k_q-1})$ also
lives in the union of $O(s^C)$ $O(2^{-\delta s/n})$-neighbourhoods of
compact hypersurfaces.  The desired cardinality
bound on the possible $J_{k_q}$ then follows from \eqref{equi} and a covering
argument.
\end{proof} 

From this proposition, we may estimate the left-hand side of
\eqref{non-deg-targ} as
$$
\sum_{(J_i)_{i \neq k_q} \in \J^{m-1}}
\int 2^{-\eps s} 2^{Ns} \prod_{i \neq k_q} \psi^+_{J_i}(x_{i-1})\ dx_0.$$
Choose an $i_0 \in \{1, \ldots, n\}$ not equal to 
$k_q$.  By applying \eqref{infty-count}
to all $J_i$ other than $J_{i_0}$, we estimate this by
$$ 
 2^{-\eps s} 2^{N(m-1)s}
\sum_{J_{i_0} \in \J} \int \psi^+_{J_{i_0}}(x_{i_0-1})\ dx_0.$$
This estimate
\eqref{non-deg-targ} then follows from \eqref{size} and \eqref{diffeo}.
This concludes the proof of \eqref{non-cancel-l1}.

\section{Proof of Theorem \ref{weak-11} continued.  The non-degenerate portion
of the integral.}

It remains to show \eqref{cancel-l1}.  By \eqref{ap-bound}, it suffices to show that
\bas
\sum_{J \in \J^m_k} 
|\int\int\int\int &G_J(x_0) F_J(x_m) \prod_{i=1}^m
\psi_{J_i}(x_{i-1}) c_{J_i}(v_i) \tilde\varphi(t_i) c_{J_i}(w_i)
\psi_{J_i}(x_i))\\
&\eta(2^{\delta s} 2^{-M} \det D_\tau(x_m)) \ dx_0
dw dt dv|
\lesssim 2^{-\eps s} \langle (T^+)^m 1, 1\rangle.
\end{align*}
By expanding out $T^+$, we see that it suffices to show that
\begin{align*}
|\int \int\int \int &G_J(x_0) F_J(x_m) \prod_{i=1}^m
\psi_{J_i}(x_{i-1}) c_{J_i}(v_i) \tilde\varphi(t_i) c_{J_i}(w_i)
\psi_{J_i}(x_i))\\
&\eta(2^{\delta s} 2^{-M} \det D_\tau(x_m)) \ dx_0
dw dt dv|\\
&\lesssim 2^{-\eps s}
\int\int\int\int \prod_{i=1}^m
(\psi^+_{J_i}(x_{i-1}) c^+_{J_i}(v_i) \varphi^+(t_i) c^+_{J_i}(w_i)
\psi^+_{J_i}(x_i))\ dx_0 dw dt dv
\end{align*}
for all $J \nearrow k$.

Fix $J \nearrow k$.
We freeze all the integration variables except for $\tau_1, \ldots,
\tau_n$, and $y$.  It thus suffices
to show that
\begin{align*}
|\int\int &G_J(x_0) \prod_{i=1}^m
\psi_{J_i}(x_{i-1}) c_{J_i}(v_i) \tilde \varphi(t_i) c_{J_i}(w_i)
\psi_{J_i}(x_i)) F_J(x_m)\\
& \eta(2^{\delta s} 2^{-M} \det D_\tau(x_m)) 
dy d\tau| \\
&\lesssim 2^{-\eps s}
\int \int \prod_{i=1}^m
\psi^+_{J_i}(x_{i-1}) c^+_{J_i}(v_i) \varphi^+(t_i) c^+_{J_i}(w_i)
\psi^+_{J_i}(x_i))\ dy d\tau
\end{align*}
uniformly in the frozen variables.

Fix all the frozen variables.  Throwing out all the factors in the above
expression which do not depend on $y$ or $\tau$, we reduce to
$$
|\int\int 
c_{J_{1}}(y)
F_J(x_m) 
a(y,\tau)\ dy d\tau| \lesssim 2^{-\eps s}
\int\int 
c^+_{J_{1}}(y) a^+(y,\tau)\ dy d\tau
$$
where
\be{a-def} a(y,\tau) = 
(\prod_{l=1}^m
\psi_{J_l}(x_{l-1}) \psi_{J_l}(x_l)) 
\eta(2^{\delta s} 2^{-M} \det D_\tau(x_m))
\prod_{q=1}^n \tilde \varphi(\tau_q)
\end{equation}
and
\be{ap-def} a^+(y,\tau) =
(\prod_{l=1}^m
\psi^+_{J_l}(x_{l-1}) \psi^+_{J_l}(x_l)) 
\prod_{q=1}^n \varphi^+(\tau_q).
\end{equation}

We now repeat the argument used to treat \eqref{cancel}.
By \eqref{div} it suffices to show that
\be{cl1-targ}
|\int\int 
\partial_{y_i} c_{J_{1}}^i(y)
F_J(x_m)  
a(y,\tau)\ 
dy d\tau| \lesssim 2^{-\eps s}
\int\int 
c^+_{J_{1}}(y)
a^+(y,\tau)
\ dy d\tau
\end{equation}
for all $i=1, \ldots, n$.

Fix $i$.  By an integration by parts, the left-hand side of
\eqref{cl1-targ} is majorized by
\be{split}
 |\int\int c_{J_{1}}^i(y)
F_J(x_m) \partial_{y_i}a(y,\tau)\ dy d\tau| +
|\int\int c_{J_{1}}^i(y)
(\partial_{y_i} F_J(x_m)) a(y,\tau)\ dy d\tau|.
\end{equation}

We now apply

\begin{lemma}  We have the pointwise estimate
\be{a-deriv}
 |(1 + 2^{\eps s} \nabla_y + \nabla_\tau) a(y,\tau)| \lesssim a^+(y,\tau).
\end{equation}
\end{lemma}

\begin{proof}
From \eqref{product}, \eqref{a-def}, \eqref{ap-def}, it suffices to verify
\bas
|(1 + 2^{\eps s} \nabla_y + \nabla_\tau) \psi_{J_l}(x_{l-1}))| &\lesssim \psi^+_{J_l}(x_{l-1})\\
|(1 + 2^{\eps s} \nabla_y + \nabla_\tau) \psi_{J_l}(x_l))| &\lesssim \psi^+_{J_l}(x_l)\\
|(1 + 2^{\eps s} \nabla_y + \nabla_\tau) \tilde \varphi(\tau_q)| &\lesssim \varphi^+(\tau_q)\\
|(1 + 2^{\eps s} \nabla_y + \nabla_\tau) \eta(2^{\delta s} 2^{-M} \det D^L_\tau(x_m))| &\lesssim 1.
\end{align*}
The first two estimates follow from \eqref{psi}, while the third is trivial.
The fourth estimate follows from the chain rule and \eqref{jac-y} providing
that $\delta \geq \eps$.
\end{proof}

From this lemma we see that the first term of \eqref{split} is acceptable.
To treat the second term, it suffices to show that
\be{split-second}
|\int (\partial_{y_i} F_J(x_m)) a(y,\tau)\ d\tau|
\lesssim 2^{-\eps s} \int a^+(y,\tau)\ d\tau
\end{equation}
uniformly in $y$.  

Fix $y$.  By Lemma \ref{chain}, we can rewrite the left-hand side as
$$
|\int \nabla_\tau F_J(x_m) \cdot
(D^L_\tau x_m)^{-1} \partial^L_{y_i} x_m a(y,\tau)\ d\tau|.
$$
Integrating by parts, we see that this is equal to
$$
|\int F_J(x_m)
\nabla_\tau \cdot ((D^L_\tau x_m)^{-1} \partial^L_{y_i} x_m a(y,\tau)) \ d\tau|.$$
Thus to show \eqref{split-second}, it suffices to verify the pointwise
estimate
$$
\|(1 + \nabla_\tau) ((D^L_\tau x_m)^{-1} \partial^L_{y_i} x_m a(y,\tau))\|
\lesssim 2^{-\eps s} a^+(y,\tau).$$
We may of course assume that $(y,\tau)$ is in the support of $a$, so that
\be{ndeg}
|\det D^L_\tau x_m| \gtrsim 2^{-\delta s} 2^M.
\end{equation}
By \eqref{a-deriv} and \eqref{product}, the left-hand side is majorized by
$$ a^+(y,\tau) |(1 + \nabla_\tau) ((D^L_\tau x_m)^{-1} \partial^L_{y_i} x_m)|$$
and so it suffices to show that
$$ \|(1 + \nabla_\tau) ((D^L_\tau x_m)^{-1} \partial^L_{y_i} x_m)\|
\lesssim 2^{-\eps s}.$$
By Cramer's rule, it suffices to show that
$$ \|(1 + \nabla_\tau) \frac{\partial^L_{\tau_1} x_m \wedge \ldots \wedge
\partial^L_{y_i} x_m \wedge \ldots \partial^L_{\tau_n} x_m}
{\det D^L_\tau x_m}\| \lesssim 2^{-\eps s}$$
for all $q$, where the numerator is the wedge product of all the
$\partial^L_{\tau_{q'}} x_m$, $q' = 1, \ldots, n$, but with the $q^{th}$ term
$\partial^L_{\tau_q} x_m$ replaced by $\partial^L_{y_i} x_m$.

Fix $q$.  From the quotient rule, \eqref{jac-y} and \eqref{ndeg} it suffices 
(if $\eps$ and $\delta$ are sufficiently small) to show that
$$ \|(1 + \nabla_\tau) (\partial^L_{\tau_1} x_m \wedge \ldots \wedge
\partial^L_{y_i} x_m \wedge \ldots \partial^L_{\tau_n} x_m)\| \lesssim 2^{-cs} 2^M$$
for some constant $c > 0$.  On the other hand, from \eqref{xl-yi} and
the inequality $j_1 \leq j_{k_q}$ arising from the hypothesis $J \nearrow k$,
we see that
$$ \| (1 + \nabla_\tau) \partial^L_{y_i} x_m \| \precsim
2^{-cs} (2^{j_{k_q} + s} \circ \one).$$
Meanwhile, from \eqref{xl-tq-v} we have
$$ \| (1 + \nabla_\tau) \partial^L_{\tau_{q'}} x_m \| \precsim
(2^{j_{k_{q'}} + s} \circ \one).$$
The desired estimate thus follows from \eqref{wedgehog} and \eqref{m-size}.
This concludes the proof of \eqref{cancel-l1} and thus of Theorem \ref{weak-11}.
\endprf

\end{document}